\documentclass[leqno]{article}
\usepackage{amssymb,latexsym,amscd,pst-all}
\newcommand{\se}[1]{{\section{#1}} {\setcounter{equation}{0}}}
\newtheorem{theorem}{Theorem}[section]
\newtheorem{lm}{Lemma}[section]
\newtheorem{prop}{Proposition}[section]

\def\k{{K\"{a}hler }}

\def\cy{{Calabi-Yau }}

\input epsf
\begin{document}
\hbadness=10000
\title{{\bf Newton polygon and string diagram}}
\author{Wei-Dong Ruan\\
Department of Mathematics\\
University of Ilinois at Chicago\\
Chicago, IL 60607\\
ruan@math.uic.edu\\}
\footnotetext{Partially supported by NSF Grant DMS-9703870 and DMS-0104150.}
\maketitle
\tableofcontents
\se{Introduction}
In this paper, we study the moment map image of algebraic curves in toric surfaces. We are particularly interested in the situations that we are able to perturb the moment map so that the moment map image of the algebraic curve is a graph. To put our problem into proper context, let's start with ${\mathbb{CP}^{2}}$.\\

Consider the natural real $n$-torus ($T^n$) action on ${\mathbb{CP}^{n}}$ given by

\[
e^{i\theta}(x) = (e^{i\theta_1}x_1, e^{i\theta_2}x_2, \cdots, e^{i\theta_n}x_n).
\]

The $T^n$ acts as symplectomorphisms with respect to the Fubini-Study \k form

\[
\omega_{\rm FS} = \partial \bar{\partial} \log(1+|x|^2).
\]

The corresponding moment map is

\[
F(x)= \left( \frac{|x_1|^2}{1+|x|^2}, \frac{|x_2|^2}{1+|x|^2}, \cdots, \frac{|x_n|^2}{1+|x|^2}\right),
\]

which is easy to see if we write $\omega_{\rm FS}$ in polar coordinates.

\[
\omega_{\rm FS} = \partial \bar{\partial} \log(1+|x|^2)=i\sum_{k=1}^n d\theta_k \wedge d\left(\frac{|x_k|^2}{1+|x|^2}\right).
\]

Notice that the moment map $F$ is a Lagrangian torus fibration and the image of the moment map $\Delta={\rm Image}(F)$ is an n-simplex.\\

In the case of ${\mathbb{CP}^2}$, $\Delta={\rm Image}(F)$ is a 2-simplex, i.e., a triangle. Let $p(z)$ be a homogeneous polynomial. $p$ defines an algebraic curve $C_p$ in ${\mathbb{CP}^2}$. We want to understand the image of $C_p$ in $\Delta$ under the moment map $F$.\\

In quantum mechanics, particle interactions are characterized by Feynman diagrams (1-dimensional graphs with some external legs). In string theory, point particles are replaced by circles (string!) and Feynman diagrams are replaced by string diagrams (Riemann surfaces with some marked points). Feynman diagrams in string theory are considered as some low energy limit of string diagrams. {\it Fattening} the Feynman diagrams by replacing points with small circles, we get the corresponding string diagrams. On the other hand, string diagrams can get ``thin" in many ways to degenerate to different Feynman diagrams.\\

Our situation is a very good analog of this picture. The complex curve $C_p$ in ${\mathbb{CP}^2}$ can be seen as a string diagram with the intersection points with the three distinguished coordinate ${\mathbb{CP}^1}$'s (that are mapped to $\partial \Delta$) as marked points. The image of $C_p$ under $F$ can be thought of as some ``fattening" of a Feynman diagram $\Gamma$ in $\Delta$ with external points in $\partial \Delta$.\\

When $p$ is of degree $d$, the genus of $C_p$ is\\
\[
g=\frac{(d-1)(d-2)}{2}.
\]
Generically, $C_p$ will intersect with ${\mathbb{CP}^1}$ at $d$ points. Ideally, the image of $C_p$ under the moment map will have $g$ holes in $\Delta$ and $d$ external points in each edge of $\Delta$. In general $F(C_p)$ can have smaller number of holes. In fact, $F(C_p)$ has at most $g$ holes. (For more detail, please see the ``Note on the literature" in the end of the introduction.)\\

In this paper we will be interested in constructing examples of $C_p$ such that $F(C_p)$ will have exactly $g$ holes in $\Delta$ and $d$ external points in each edge of $\Delta$. Namely, the case when $F(C_p)$ resembles classical Feynman diagrams the most. (Sort of the most classical string diagram.) These examples will be constructed for any degree in section 2.\\

Our interest on this problem comes from our work on Lagrangian torus fibration of \cy manifold and mirror symmetry. In \cite{lag1,lag2,lag3}, we mainly concern the case of quintic curves in ${\mathbb{CP}^2}$. The generalization to curves in toric surfaces will be useful in \cite{tor,ci}. The algebraic curves and their images under the moment map arise as the singular set and singular locus of our Lagrangian torus fibrations.\\

As we mentioned, $F(C_p)$ can be rather chaotic for general curve $C_p$. The condition for $F(C_p)$ to resemble a classical Feynman diagram is related to the concept of ``\textsf{near the large complex limit}'', which is explained in section 3. (Through discussion with Qin Jing, it is apparent that ``near the large complex limit'' is equivalent to near the 0-dimensional strata in $\overline{{\cal M}_g}$, the moduli space of stable curves of genus $g$. These points in $\overline{{\cal M}_g}$ are represented by stable curves, whose irreducible components are all ${\mathbb{CP}^1}$ with three marked points.) It turns out that our construction of ``graph like'' string diagrams for curves in ${\mathbb{CP}^2}$ can be generalized to curves in general 2-dimensional toric varieties using localization technique. More precisely, \textsf{in the moduli space of curves in a general 2-dimensional toric variety, when the curve $C_p$ is close enough to the so-called ``large complex limit''(analogous to classical limit in physics) in suitable sense,  $F(C_p)$ will resemble a fattening of a classical Feynman diagram.} This result will be made precise and proved in theorems \ref{ca} and \ref{cg} of section 3. Examples constructed in section 2 are special cases of this general construction.\\

One advantage of string theory over classical quantum mechanic is that the string diagrams (marked Riemann surfaces) are more natural than Feynman diagrams (graphs). For instance, one particular topological type of string diagram under different classical limit can degenerate into very different Feynman diagrams, therefore unifying them. In our construction, there is a natural partition of the moduli space of curves such that in different part the limiting Feynman diagrams are different. We will discuss this natural partition of the moduli space and different limiting Feynman diagrams also in section 3.\\

Of course, ideally, it will be interesting if $F(C_p)$ is actually a 1-dimensional Feynman diagram $\Gamma$ in $\Delta$. This will not be true for the moment map $F$. A natural question is: ``Can one perturb the moment map $F$ to $\hat{F}$ so that $\hat{F}(C_p) = \Gamma$?" (Notice that the moment map of a torus action is equivalent to a Lagrangian torus fibration. We will use the two concepts interchangeably in this paper.) Such perturbation is not possible in the smooth category. But \textsf{when $F(C_p)$ resembles a classical Feynman diagram $\Gamma$ close enough, we can perturb $F$ suitably as a moment map, so that the perturbed moment map $\hat{F}$ is piecewise smooth and satisfies $\hat{F}(C_p) = \Gamma$.} This perturbation construction is explicitly done for the case of line in ${\mathbb{CP}^2}$ in section 4 (theorems \ref{dc} and \ref{dt}). The general case is dealt with in section 5 (theorems \ref{dd} and \ref{ee}) combining the localization technique in section 3 and the perturbation technique in section 4. (In particular, optimal smoothness for $\hat{F}$ is achieved in theorems \ref{dt} and \ref{ee}.)\\

{\bf Note on the literature:} Our work on Newton polygon and string diagram was motivated by and was an important ingredient of our construction of Lagrangian torus fibrations of Calabi-Yau manifolds (\cite{lag1,lag2,lag3,tor,ci}). After reading my preprint, Prof. Y.-G. Oh pointed out to me the work of G. Mikhalkin (\cite{M}), through which I was able to find the literature of our problem. The image of curves under the moment map was first investigated in \cite{G}, where it was called ``amoeba". Legs of ``amoeba" are already understood in \cite{G}. The problem of determining holes in ``amoeba" was posed in Remark 1.10 of page 198 in \cite{G} as a difficult and interesting problem. Work of G. Mikhalkin (\cite{M}) that was published in 2000 and works (\cite{F}, \cite{V}, etc.) mentioned in its reference point out some previous progress on this problem of determining holes in ``amoeba" aimed at very different applications, which nevertheless is very closely related to our work. Most of the ideas in section 2 and 3 are not new and appeared in one form or the other in these previous works mentioned. For example, our localization technique used in section 3 closely resemble the curve patching idea of Viro (which apparently appeared much earlier) in different context as described in \cite{M}. Due to different purposes, our approach and results are of somewhat distinctive flavor. To our knowledge, our discussion in section 4 and 5 on symplectic deformation to Lagrangian fibrations with the image of curve being graph, which is important for our applications, was not discussed before and is essentially new. I also want to mention that according to the description in \cite{M} of a result of Forsberg et al. \cite{F}, one can derive that there are at most $g = \frac{(d-1)(d-2)}{2}$ holes in $F(C_p)$ for degree $d$ curve $C_p$, which I initially conjectured to be true.\\

{\bf Note on the figures:} The figures of moment map images of curves as fattening of graphs in this paper are somewhat idealized topological illustration. Some part of the edges of the image that are straight or convex could be curved or concave in more accurate picture. Of course, such inaccuracy will not affect our mathematical argument and the fact that moment map images of curves are fattening of graphs.\\

{\bf Notion of convexity:} A function $y = f(x)$ will be called convex if the set $\{(x,y): y \geq f(x)\}$ is convex. We are aware such functions have been called concave by some authors.

\se{The construction for curves in $\mathbb{CP}^2$}
To understand our problem better, let us look at the example of Fermat type polynomial\\
\[
p=z_1^d + z_2^d + z_3^d.
\]
It is not hard to see that for any $d$, $F(C_p)$ will look like a curved triangle with only one external point in each edge of $\Delta$ and no hole at all. (This example is in a sense a string diagram with the most quantum effect.)\\
\begin{center}
\begin{picture}(200,50)(0,130)
\thicklines
\put(64,158){\line(3,-1){12}}
\put(76,154){\line(1,0){12}}
\put(100,158){\line(-3,-1){12}}

\put(18,30)
{\put(46,128){\line(4,-3){10}}
\put(56,120.5){\line(3,-5){6}}
\put(64,98){\line(-1,6){2.1}}}

\put(146,30)
{\put(-46,128){\line(-4,-3){10}}
\put(-56,120.5){\line(-3,-5){6}}
\put(-64,98){\line(1,6){2.1}}}

\thinlines
\multiput(82,149)(36,0){1}{\line(2,1){18}}
\multiput(82,149)(36,0){1}{\line(-2,1){18}}
\multiput(82,128)(36,0){1}{\line(0,1){21}}
\put(46,128){\line(1,0){72}}
\put(46,128){\line(3,5){36}}
\put(118,128){\line(-3,5){36}}
\end{picture}
\end{center}
\begin{center}
\stepcounter{figure}
Figure \thefigure: $F(C_p)$ of $p(z) =z_1^d + z_2^d + z_3^d$.
\end{center}
From this example, it is not hard to imagine that for most polynomials, chances are the number of holes will be much less than $g$. Any attempt to construct examples with the maximal number of holes will need special care, especially if one wants the construction for general degree $d$.\\

Let $[z]=[z_1,z_2,z_3]$ be the homogeneous coordinate of ${\mathbb{CP}^2}$. Then a general homogeneous polynomial of degree $d$ in $z$ can be expressed as\\
\[
p(z) = \sum_{I\in N^d}a_I z^I,
\]
where\\
\[
N^d = \{I=(i_1,i_2,i_3)\in {\mathbb{Z}^3}| |I|=i_1+i_2+i_3=d, I\geq 0\}
\]
is the Newton polygon of degree $d$ homogeneous polynomials. In our case $N^d$ is a triangle with $d+1$ lattice points on each side. Denote $E=(1,1,1)$.\\

To describe our construction, let us first notice that $N^d$ can be naturally decomposed as a union of "hollow" triangles as follows:\\
\[
N^d = \bigcup_{k=0}^{[d/3]} N^d_k,
\]
where\\
\[
N^d_k = \{I\in N^d | I\geq kE,\ I\not\geq (k+1)E\}.
\]
On the other hand, the map $I\rightarrow I+E$ naturally defines an embedding $i: N^d \rightarrow N^{d+3}$. From this point of view, $N^d_0 = N^d\backslash N^{d-3}$ and $N^d_k = i(N^{d-3}_{k-1})$ for $k\geq 1$.\\

When $d=1$, $g=0$ and a generic degree 1 polynomial can be reduced to\\
\[
p=z_1 + z_2 +z_3.
\]
$F(C_p)$ is a triangle with vertices as middle points of edges of $\Delta$. This clearly satisfies our requirement, namely, with $g=0$ holes.\\

For $d\geq 2$, the first problem is to make sure that the external points are distinct and as far apart as possible. For this purpose, we want to consider homogeneous polynomials with two variables. A nice design is to consider\\
\[
q_d(z_1,z_2) = \prod_{i=1}^d(z_1+ t_iz_2) = \sum_{i=0}^d b_i z_1^{d-i}z_2^i
\]
such that $t_{d-i+1}=\frac{1}{t_i}\geq 1$. Then $b_0=b_d=1$ and $b_{d-i} = b_i \geq 1$ for $i \geq 1$. We can adjust $t_i$ for $1\leq i \leq [d/2]$ suitably to make them far apart. (For example, one may assume $F(t_i)= \frac{2i-1}{2d}$, where $F(t)=\frac{t^2}{1+t^2}$ is the moment map for $n=1$.) Now we can define a degree $d$ homogeneous polynomial in three variables such that coefficients along each edge of $N_d$ is assigned according to $q$ and coefficients in the interior of $N_d$ vanish. We will still denote this polynomial by $q_d$. Then we have\\
\begin{theorem}
\label{ba}
For
\[
p_d(z) = \sum_{k=0}^{[d/3]} c_k q_{d-3k}(z) z^{kE},
\]
where $c_0=1$ and $c_k>0$, if $c_k$ is big enough compared to $c_{k-1}$, then $F(C_{p_d})$ has exactly $g$ holes and $d$ external points in each edge of $\Delta$.
\end{theorem}
Before proving the theorem, let us analyze some examples that give us better understanding of the theorem. The following are examples of Feynman diagrams corresponding to $1\leq d \leq 5$. (The case $d=5$ is the case we are interested in mirror symmetry.) The image of the corresponding string diagrams under the moment map $F$ are some fattened version of these Feynman diagram. For example one can see genus of the corresponding Riemann surfaces from these diagrams. When $d=1,2$ there are no holes in the diagram and genus equal to zero. When $d=5$ there are $6$ holes and the corresponding Riemann surface are genus $6$ curves. These diagrams give a very nice interpretation of genus formula for planar curve. (In my opinion, also a good way to remember it!)\\

To justify our claim, we first analyze it case by case. When $d=1$, $g=0$ and a generic polynomial can be reduced to\\
\[
p_1(z)=z_1 + z_2 +z_3.
\]
$F(C_p)$ is a triangle with vertices as middle points of edges of $\Delta$. This clearly is a fattened version of the first diagram in the next picture.\\

When $d=2$, $g=0$ and we may take polynomial\\
\[
p_2(z) = (z_1^2 + z_2^2 +z_3^2) + \frac{5}{2}(z_1z_2 + z_2z_3 +z_3z_1).
\]
When $z_3=0$\\
\[
q_2(z_1,z_2) = z_1^2 + z_2^2 + \frac{5}{2}z_1z_2 = (z_1 + \frac{1}{2}z_2)(z_1 + 2z_2).
\]
Image of $\{q_2(z_1,z_2)=0\}$ under $F$ are two lines coming out of the edge $r_3=0$ starting from the two points $\frac{r_1}{r_2} = 2, \frac{1}{2}$. When $z_3$ is small, $p_2(z)$ is a small perturbation of $q_2(z_1,z_2)$. By this argument, it is clear that $F(C_2)$ is a fattening of the second diagram in the following picture near the boundary of the triangle. Since in our case $g=0$. It is not hard to conceive or (if you are more strict) to find a way to prove that $F(C_2)$ is a fattening of the second diagram in the following picture.\\
\begin{center}
\begin{picture}(200,80)(-90,90)
\put(-150,-10){\thicklines
\put(64,158){\line(3,-1){12}}
\put(76,154){\line(1,0){12}}
\put(100,158){\line(-3,-1){12}}

\put(18,30)
{\put(46,128){\line(4,-3){10}}
\put(56,120.5){\line(3,-5){6}}
\put(64,98){\line(-1,6){2.1}}}

\put(146,30)
{\put(-46,128){\line(-4,-3){10}}
\put(-56,120.5){\line(-3,-5){6}}
\put(-64,98){\line(1,6){2.1}}}

\thinlines
\multiput(82,149)(36,0){1}{\line(2,1){18}}
\multiput(82,149)(36,0){1}{\line(-2,1){18}}
\multiput(82,128)(36,0){1}{\line(0,1){21}}
\put(46,128){\line(1,0){72}}
\put(46,128){\line(3,5){36}}
\put(118,128){\line(-3,5){36}}}

\thicklines
\put(64,158){\line(3,-1){12}}
\put(76,154){\line(1,0){12}}
\put(100,158){\line(-3,-1){12}}

\put(64,158){\line(4,-3){15.2}}
\put(79.2,146.6){\line(0,-1){17}}
\put(46,128){\line(3,-1){18}}
\put(64,122){\line(2,1){15.4}}

\put(164,0)
{\put(-64,158){\line(-4,-3){15.2}}
\put(-79.2,146.6){\line(0,-1){17}}
\put(-46,128){\line(-3,-1){18}}
\put(-64,122){\line(-2,1){15.4}}}

\put(46,128){\line(4,-3){10}}
\put(56,120.5){\line(3,-5){6}}
\put(64,98){\line(-1,6){2.1}}

\put(164,0)
{\put(-46,128){\line(-4,-3){10}}
\put(-56,120.5){\line(-3,-5){6}}
\put(-64,98){\line(1,6){2.1}}}

\put(64,98){\line(1,6){3.3}}
\put(100,98){\line(-1,6){3.3}}
\put(82,125){\line(2,-1){15}}
\put(82,125){\line(-2,-1){15}}

\thinlines
\multiput(82,149)(36,0){1}{\line(2,1){18}}
\multiput(82,149)(36,0){1}{\line(-2,1){18}}
\multiput(82,128)(36,0){1}{\line(0,1){21}}
\multiput(46,128)(36,0){2}{\line(2,-1){18}}
\multiput(82,128)(36,0){2}{\line(-2,-1){18}}
\multiput(64,98)(36,0){2}{\line(0,1){21}}
\put(28,98){\line(1,0){108}}
\put(28,98){\line(3,5){54}}
\put(136,98){\line(-3,5){54}}
\end{picture}
\end{center}
\begin{center}
\refstepcounter{figure} \label{fig1}
Figure \thefigure: degree $d=1,2$
\end{center}
When $d=3$, $g=1$. We can consider\\
\[
p_3(z) = (z_1^3 + z_2^3 +z_3^3) + \frac{7}{2}(z_1^2z_2 + z_2^2z_3 +z_3^2z_1 + z_1z_2^2 + z_2z_3^2 +z_3z_1^2) + bz_1z_2z_3.
\]
We can use similar idea as in the previous case to explain the behavior of $F(C_3)$ near the edges. The main point for this case is to explain how the hole in the center arises. For this purpose, we introduce the following function\\
\[
\rho_p(r)=\inf_{F(z)=r}|p(z)|.
\]
This function takes non-negative value, and\\
\[
F(C_p)=\{[r]|\rho_p(r)=0\}.
\]
$\rho$ also satisfies\\
\[
\rho_{p_1p_2} = \rho_{p_1} \rho_{p_2}
\]
\[
\rho_{p_1+p_2} \leq \rho_{p_1} + \rho_{p_2}
\]
An important thing to notice is that $\rho_{z_1z_2z_3}(r) = r_1r_2r_3$ is a function that vanishes at  the edges of the triangle and not vanishing anywhere in the interior of the triangle, sort of a bump function. When $b$ is large, $\rho_{p_3}$ will be dominated by $br_1r_2r_3$ away from the edges, which will be positive around center. Therefore $F(C_3)$ will have a hole in the center, which becomes large when $b$ gets large.\\
\begin{center}
\begin{picture}(200,120)(-20,60)
\thicklines
\put(64,158){\line(3,-1){12}}
\put(76,154){\line(1,0){12}}
\put(100,158){\line(-3,-1){12}}

\multiput(0,0)(-18,-30){2}
{\put(64,158){\line(4,-3){15.2}}
\put(79.2,146.6){\line(0,-1){17}}
\put(46,128){\line(3,-1){18}}
\put(64,122){\line(2,1){15.4}}}

\multiput(164,0)(18,-30){2}
{\put(-64,158){\line(-4,-3){15.2}}
\put(-79.2,146.6){\line(0,-1){17}}
\put(-46,128){\line(-3,-1){18}}
\put(-64,122){\line(-2,1){15.4}}}

\put(-18,-30)
{\put(46,128){\line(4,-3){10}}
\put(56,120.5){\line(3,-5){6}}
\put(64,98){\line(-1,6){2.1}}}

\put(182,-30)
{\put(-46,128){\line(-4,-3){10}}
\put(-56,120.5){\line(-3,-5){6}}
\put(-64,98){\line(1,6){2.1}}}

\multiput(-18,-30)(36,0){2}
{\put(64,98){\line(1,6){3.3}}
\put(100,98){\line(-1,6){3.3}}
\put(82,125){\line(2,-1){15}}
\put(82,125){\line(-2,-1){15}}}

\put(82,108.5){\circle{30}}

\thinlines
\multiput(82,149)(36,0){1}{\line(2,1){18}}
\multiput(82,149)(36,0){1}{\line(-2,1){18}}
\multiput(82,128)(36,0){1}{\line(0,1){21}}
\multiput(46,128)(36,0){2}{\line(2,-1){18}}
\multiput(82,128)(36,0){2}{\line(-2,-1){18}}
\multiput(64,98)(36,0){2}{\line(0,1){21}}
\multiput(28,98)(36,0){3}{\line(2,-1){18}}
\multiput(64,98)(36,0){3}{\line(-2,-1){18}}
\multiput(46,68)(36,0){3}{\line(0,1){21}}
\put(10,68){\line(1,0){144}}
\put(10,68){\line(3,5){72}}
\put(154,68){\line(-3,5){72}}
\end{picture}
\end{center}
\begin{center}
\stepcounter{figure}
Figure \thefigure: degree $d=3$
\end{center}
When $d=4$, $g=3$. We can consider\\
\[
p_4(z) = q_4(z) + bz_1z_2z_3p_1(z).
\]
The key point is to understand how the three holes appear. For this purpose, we need to go back to the case when $d=1$. Notice that $\rho_{z_1z_2z_3p_1}=r_1r_2r_3\rho_{p_1}$ is positive in the three regions as indicated in the diagram for $d=1$, and it is zero at the boundary of the three regions. When $b$ is large, this term dominates $\rho_{p_4}$ in the interior of the triangle and produces the three holes. Similar discussion as before implies that $q_4$ will take care of edges.\\
\begin{center}
\begin{picture}(200,150)(-30,30)
\thicklines
\put(64,158){\line(3,-1){12}}
\put(76,154){\line(1,0){12}}
\put(100,158){\line(-3,-1){12}}

\multiput(0,0)(-18,-30){3}
{\put(64,158){\line(4,-3){15.2}}
\put(79.2,146.6){\line(0,-1){17}}
\put(46,128){\line(3,-1){18}}
\put(64,122){\line(2,1){15.4}}}

\multiput(164,0)(18,-30){3}
{\put(-64,158){\line(-4,-3){15.2}}
\put(-79.2,146.6){\line(0,-1){17}}
\put(-46,128){\line(-3,-1){18}}
\put(-64,122){\line(-2,1){15.4}}}

\put(-36,-60)
{\put(46,128){\line(4,-3){10}}
\put(56,120.5){\line(3,-5){6}}
\put(64,98){\line(-1,6){2.1}}}

\put(200,-60)
{\put(-46,128){\line(-4,-3){10}}
\put(-56,120.5){\line(-3,-5){6}}
\put(-64,98){\line(1,6){2.1}}}

\multiput(-36,-60)(36,0){3}
{\put(64,98){\line(1,6){3.3}}
\put(100,98){\line(-1,6){3.3}}
\put(82,125){\line(2,-1){15}}
\put(82,125){\line(-2,-1){15}}}

\multiput(82,108.5)(18,-30){2}{\circle{30}}
\multiput(64,78.5)(18,-30){1}{\circle{30}}

\thinlines
\multiput(82,149)(36,0){1}{\line(2,1){18}}
\multiput(82,149)(36,0){1}{\line(-2,1){18}}
\multiput(82,128)(36,0){1}{\line(0,1){21}}
\multiput(46,128)(36,0){2}{\line(2,-1){18}}
\multiput(82,128)(36,0){2}{\line(-2,-1){18}}
\multiput(64,98)(36,0){2}{\line(0,1){21}}
\multiput(28,98)(36,0){3}{\line(2,-1){18}}
\multiput(64,98)(36,0){3}{\line(-2,-1){18}}
\multiput(46,68)(36,0){3}{\line(0,1){21}}
\multiput(10,68)(36,0){4}{\line(2,-1){18}}
\multiput(46,68)(36,0){4}{\line(-2,-1){18}}
\multiput(28,38)(36,0){4}{\line(0,1){21}}
\put(-8,38){\line(1,0){180}}
\put(-8,38){\line(3,5){90}}
\put(172,38){\line(-3,5){90}}
\end{picture}
\end{center}
\begin{center}
\stepcounter{figure}
Figure \thefigure: degree $d=4$
\end{center}
When $d=5$, $g=6$. We need to go back to the case $d=2$. The discussion is very similar to the previous case, we will omit.\\
\begin{center}
\begin{picture}(200,180)(-30,0)
\thicklines
\put(64,158){\line(3,-1){12}}
\put(76,154){\line(1,0){12}}
\put(100,158){\line(-3,-1){12}}

\multiput(0,0)(-18,-30){4}
{\put(64,158){\line(4,-3){15.2}}
\put(79.2,146.6){\line(0,-1){17}}
\put(46,128){\line(3,-1){18}}
\put(64,122){\line(2,1){15.4}}}

\multiput(164,0)(18,-30){4}
{\put(-64,158){\line(-4,-3){15.2}}
\put(-79.2,146.6){\line(0,-1){17}}
\put(-46,128){\line(-3,-1){18}}
\put(-64,122){\line(-2,1){15.4}}}

\put(-54,-90)
{\put(46,128){\line(4,-3){10}}
\put(56,120.5){\line(3,-5){6}}
\put(64,98){\line(-1,6){2.1}}}

\put(218,-90)
{\put(-46,128){\line(-4,-3){10}}
\put(-56,120.5){\line(-3,-5){6}}
\put(-64,98){\line(1,6){2.1}}}

\multiput(-54,-90)(36,0){4}
{\put(64,98){\line(1,6){3.3}}
\put(100,98){\line(-1,6){3.3}}
\put(82,125){\line(2,-1){15}}
\put(82,125){\line(-2,-1){15}}}

\multiput(82,108.5)(18,-30){3}{\circle{30}}
\multiput(64,78.5)(18,-30){2}{\circle{30}}
\multiput(46,48.5)(18,-30){1}{\circle{30}}

\thinlines
\multiput(82,149)(36,0){1}{\line(2,1){18}}
\multiput(82,149)(36,0){1}{\line(-2,1){18}}
\multiput(82,128)(36,0){1}{\line(0,1){21}}
\multiput(46,128)(36,0){2}{\line(2,-1){18}}
\multiput(82,128)(36,0){2}{\line(-2,-1){18}}
\multiput(64,98)(36,0){2}{\line(0,1){21}}
\multiput(28,98)(36,0){3}{\line(2,-1){18}}
\multiput(64,98)(36,0){3}{\line(-2,-1){18}}
\multiput(46,68)(36,0){3}{\line(0,1){21}}
\multiput(10,68)(36,0){4}{\line(2,-1){18}}
\multiput(46,68)(36,0){4}{\line(-2,-1){18}}
\multiput(28,38)(36,0){4}{\line(0,1){21}}
\multiput(-8,38)(36,0){5}{\line(2,-1){18}}
\multiput(28,38)(36,0){5}{\line(-2,-1){18}}
\multiput(10,8)(36,0){5}{\line(0,1){21}}
\put(-26,8){\line(1,0){216}}
\put(-26,8){\line(3,5){108}}
\put(190,8){\line(-3,5){108}}
\end{picture}
\end{center}
\begin{center}
\refstepcounter{figure} \label{fig2}
Figure \thefigure: degree $d=5$
\end{center}
{\bf Proof of theorem \ref{ba}:} We prove by induction. For this purpose, notice that we can define $p_d(z)$ alternatively by induction\\
\[
p_d(z)=q_d(z) + b_d z^Ep_{d-3}(z).
\]
We need to show that when $b_d$ are large enough for any $d$, $F(C_{p_d})$ will have $g=\frac{(d-1)(d-2)}{2}$ holes and $d$ external legs in each edge.\\

Assume above statement is true for $p_{d-3}(z)$, then $F(C_{p_{d-3}})$ will have $g=\frac{(d-4)(d-5)}{2}$ holes and $d-3$ external legs in each edge. It is easy to see that $F(C_{z^Ep_{d-3}})$ will have $g=\frac{(d-4)(d-5)}{2}$ interior holes and $3(d-3)$ side holes that are partly bounded by edges. We are expecting that by adding $q_d$ term, side holes will become interior hole and there will be $d$ external leges on each edge.\\

Discussion in previous special examples will more or less do this. Here we can do better. We can actually write down explicitly the behavior of $F(C_p)$ near edges. For example, near the edge $z_3=0$, $p_d(z)=0$ can be rewritten as\\
\[
z_3= -\frac{q_d(z)}{b_dz_1z_2p_{d-3}(z)}.
\]
This is a graph over the coordinate line $z_3=0$ within say $|z_3|\leq \epsilon\min(|z_1|,|z_2|)$ and away from $z_1=0$, $z_2=0$ and $d-3$ leg points of $p_{d-3}(z)$. It will be clearer to discuss under local coordinate say $x_1 = \frac{z_1}{z_2}$, $x_3 = \frac{z_3}{z_2}$. We will use the same symbol for homogeneous polynomials and the corresponding inhomogeneous polynomials. Then under this inhomogeneous coordinate\\
\[
x_3= -\frac{q_d(x_1,x_3)}{b_dx_1p_{d-3}(x_1,x_3)}.
\]
Asymptotically, near $x_3=0$\\
\[
x_3= -\frac{q_d(x_1)}{b_dx_1q_{d-3}(x_1)}.
\]
From previous notation $q_d(x_1)=q_d(x_1,0) = p_d(x_1,0)$, and\\
\[
q_d(x_1) = \prod_{i=1}^d(x_1- t_{d,i}).
\]
Recall that we require $|t_{d,i}|$ to be as far apart as possible for different $d,i$. From this explicit expression, it is easy to see that near $z_3=0$ (say $|x_3|\leq\epsilon$) and away from $z_2=0$, $C_{p_d}$ is a graph over the ${\mathbb{CP}^1}$ ($z_3=0$) away from disks\\
\[
|x_1 - t_{d-3,i}|\leq \frac{q_d(t_{d-3,i})}{t_{d-3,i}q'_{d-3}(t_{d-3,i})}\frac{1}{\epsilon b_d} \ \  {\rm for}\ 1\leq i \leq d-3,
\]
and
\[
|x_1|\leq \frac{q_d(0)}{q_{d-3}(0)}\frac{1}{\epsilon b_d}.
\]
Recall $b_d$ is supposed to be large. Here we further require the choice of $\epsilon$ to satisfy$\epsilon$ is small and $\epsilon b_d$ is large. Therefore, all these holes are very small. It is easy to see that the $d-3$ small circles centered around the roots of $q_{d-3}$ will connect with $d-3$ legs of $C_{p_{d-3}}$. In this way, the side holes of $F(C_{z^Ep_{d-3}})$ will become interior holes of $F(C_p)$. Together with original interior holes they add up to\\
\[
\frac{(d-4)(d-5)}{2}+3(d-3) = \frac{(d-1)(d-2)}{2} = g_d
\]
interior holes for $F(C_{p_d})$. $d$ zeros of $q_d$ along each edge will produce for us the $d$ exterior legs on each edge. Namely $F(C_d)$ is fattening of the Feynman diagrams as described in previous pictures.
\hfill$\Box$\\

\se{Newton polygon and string diagram}
The result in the previous section is actually special cases of a more general result on curves in toric surfaces. When the coefficients of the defining equation of a curve in a general toric surface satisfy certain convexity conditions (in physical term: ``near the large complex limit"), the moment map image (amoeba) of the curve in the toric surface will also resemble the fattening of a graph. The key idea that enables such generalization is the so-called ``localization technique" that reduces the amoeba of our curve near the large complex limit locally to the amoeba of a line, which is well understood.\\

We start with toric terminologies. Let $M$ be a rank 2 lattice and $N = M^\vee$ denotes the dual lattice. For any $\mathbb{Z}$-module $\mathbb{A}$, let $N_{\mathbb{A}} = N \otimes_\mathbb{Z} \mathbb{A}$. Given an integral polygon $\Delta\subset M$, we can naturally associate a fan $\Sigma$ by the construction of normal cones. For a face $\alpha$ of the polygon $\Delta$, define the normal cone of $\alpha$\\
\[
\sigma_\alpha :=\{n\in N|\langle m',n\rangle \leq \langle m,n\rangle \ \ {\rm for\ all\ } m'\in \alpha,\ m\in \Delta \}.
\]\\
Let $\Sigma$ denote the fan that consists of all these normal cones. We are interested in the corresponding toric variety $P_\Sigma$. Let $\Sigma(1)$ denote the collection of one dimensional cones in the fan $\Sigma$, then any $\sigma\in\Sigma(1)$ determines a $N_{\mathbb{C}^*}$-invariant Weil divisor $D_{\sigma}$.\\

For $m \in M$, $s_m = e^{\langle m,n\rangle}$ defines a monomial function on $N_{\mathbb{C}^*}$ that extends to a meromorphic function on $P_\Sigma$. Let $e_\sigma$ denote the unique primitive element in $\sigma \in\Sigma(1)$. The Cartier divisor

\[
(s_m) = \sum_{\sigma\in \Sigma(1)} \langle m,e_\sigma\rangle D_\sigma.
\]

Consider the divisor\\
\[
D_\Delta = \sum_{\sigma\in \Sigma(1)} l_\sigma D_\sigma, \mbox{ where } l_\sigma = -\inf_{m\in \Delta} \langle m,e_\sigma\rangle.
\]\\
The corresponding line bundle $L_\Delta = {\cal O}(D_\Delta)$ can be characterized by the piecewise linear function $p_{\Delta}$ on $N$ that satisfies $p_\Delta(e_\sigma)=l_\sigma$ for any $\sigma\in \Sigma(1)$. It is easy to see that $p_\Delta$ is strongly convex with respect to the fan $\Sigma$, hence $L_\Delta$ is ample on $P_\Sigma$. Since $(s_m) + D_\Delta$ is effective if and only if $m\in \Delta$, $\{s_m\}_{m \in \Delta}$ can be identified with the set of $N_{\mathbb{C}^*}$-invariant holomorphic sections of $L_\Delta$. In this sense, the polygon $\Delta$ is usually called the {\bf Newton polygon} of the line bundle $L_\Delta$ on $P_\Sigma$. A general section of $L_\Delta$ can be expressed as\\
\[
s = \sum_{m\in \Delta} a_ms_m.
\]\\
$C_s = s^{-1}(0)$ is a curve in $P_\Sigma$. We can consider the image of the curve $C_s$ under some moment map of $P_\Sigma$. The problem we are interested in is when this image will form a fattening of a graph. The case discussed in the last section is a special case of this problem, corresponding to the situation of $P_\Sigma \cong \mathbb{CP}^2$ and $L_\Delta \cong {\cal O}(k)$.\\

With $w= \{w_m\}_{m\in \Delta} \in \tilde{N}_0 \cong \mathbb{Z}^{\Delta}$, we can define an action of $\delta\in {\mathbb{R}}_+$ on sections of $L_\Delta$.

\[
s^{\delta^w} = \delta(s)= \sum_{m\in \Delta} (\delta^{w_m}a_m)s_m.
\]
\[
A = \{\{l + \langle m,n\rangle\}_{m\in \Delta} \in \tilde{N}_0: (l,n) \in N^+ = \mathbb{Z} \oplus N\} \subset \tilde{N}_0
\]

is the sublattice of affine functions on $\Delta$. An element $[w] \in \tilde{N} = \tilde{N}_0/A$ can be viewed as an equivalent class of $\mathbb{Z}$-valued functions $w=(w_m)_{m\in \Delta}$ on $\Delta$ modulo the restriction of affine function on $M$.\\

When $w= \{w_m\}_{m\in \Delta} \in \tilde{N}_0$ is a strictly convex function on $\Delta$, $w$ determines a simplicial decomposition $Z$ of $\Delta$. Clearly every representative of $[w] \in \tilde{N}$ determines the same simplicial decomposition $Z$ of $\Delta$. Let $\tilde{S}$ (resp. $\tilde{S}^{\rm top}$) be the set of $S\subset \Delta$ that forms a simplex (resp. top dimensional simplex) containing no other integral points. Then $Z$ can be regarded as a subset of $\tilde{S}$. Let $Z^{\rm top} = Z \cap \tilde{S}^{\rm top}$.\\

From now on, assume $|a_m|=1$ for all $m\in \Delta$. $|s_m| = |e^{\langle m,n\rangle}|$ is a function on $N_{\mathbb{C}^*} \subset P_\Sigma$. Let

\[
h_{\delta^w}=\log |s^{\delta^w} |_{\Delta}^2, \mbox{ where } |s^{\delta^w} |_{\Delta}^2 = \sum_{m\in\Delta}|s^{\delta^w}|_m^2,\ |s^{\delta^w}|_m = |s^{\delta^w}_m| = |\delta^{w_m} s_m|.
\]

$\omega_{\delta^w} = \partial \bar{\partial} h_{\delta^w}$ naturally defines a $N_{\mathbb{S}}$-invariant \k form on $P_{\Sigma}$, where $\mathbb{S}$ denotes the unit circle in $\mathbb{C}^*$ as $\mathbb{Z}$-submodule.\\

Choose a basis $n_1, n_2$ of $N$, then $n\in N_\mathbb{C}$ can be expressed as

\[
n =\sum_{k=1}^2 (\log x_k)n_k =\sum_{k=1}^2 (\log r_k + i\theta_k)n_k.
\]

Under this local coordinate, the \k form $\omega_{\delta^w}$ can be expressed as

\[
\omega_{\delta^w} = \partial \bar{\partial} h_{\delta^w} = i\sum_{k=1}^2 d\theta_k \wedge d h_k, \mbox{ where } h_k=|x_k|^2\frac{\partial h_{\delta^w}}{\partial |x_k|^2}.
\]

It is straightforward to compute that

\[
h_k = \sum_{m\in\Delta} \langle m,n_k\rangle \rho_m, \mbox{ where } \rho_m = \frac{|s^{\delta^w} |_m^2}{|s^{\delta^w} |_{\Delta}^2}.
\]

Consequently,

\[
\omega_{\delta^w} = i\sum_{m\in\Delta} d\langle m,\theta\rangle \wedge d\rho_m, \mbox{ where } \theta = \sum_{k=1}^2 \theta_k n_k.
\]
\begin{lm}
The moment map is\\
\[
F_{\delta^w}(x) = \sum_{m\in\Delta} \rho_m(x) m.
\]\\
which maps $P_\Sigma$ to $\Delta$.
\hfill$\Box$\\
\end{lm}
By this map, $N_{\mathbb{S}}$-invariant functions $h$, $h_k$, $\rho_m$ on $P_\Sigma$ can all be viewed as functions on $\Delta$. We have\\
\begin{lm}
\label{cb}
$\rho_m$ as a function on $\Delta$ achieves its maximum exactly at $m\in \Delta$.
\end{lm}
{\bf Proof:} By $\displaystyle x_k\frac{\partial |s^{\delta^w} |_m^2}{\partial x_k} = \langle m,n_k\rangle |s^{\delta^w} |_m^2$, $\rho_m$ achieves maximal implies\\

\[
\sum_{k=1}^2 x_k\frac{\partial \rho_m}{\partial x_k}m_k = \rho_m \sum_{m'\in\Delta}\sum_{k=1}^2(\langle m,n_k\rangle - \langle m',n_k\rangle)\rho_{m'}m_k
\]
\[
= \rho_m \sum_{m'\in\Delta}(m - m')\rho_{m'} = \rho_m(m - F_{\delta^w}(x))=0.
\]\\
Therefore $F_{\delta^w}(x) =m$ when $\rho_m$ achieves maximal.
\hfill$\Box$\\
\begin{lm}
\label{ce}
For any subset $S \subset \Delta$, $\rho_S = \displaystyle \sum_{m\in S} \rho_m$ as a function on $\Delta$ achieves maximum in the convex hull of $S$. At the maximal point of $\rho_S$

\[
F_{\delta^w}(x) = \sum_{m\in S}\rho^S_m m = \sum_{m\not\in S}\rho^{S^c}_m m, \mbox{ where } \rho^S_m = \frac{\rho_m}{\rho_S},\ S^c = \Delta \setminus S.
\]
\end{lm}
{\bf Proof:} By $\displaystyle x_k\frac{\partial |s^{\delta^w} |_m^2}{\partial x_k} = \langle m,n_k\rangle |s^{\delta^w} |_m^2$, $\rho_S = \displaystyle \sum_{m\in S} \rho_m$ achieves maximal implies\\

\[
\sum_{m\in S}\sum_{k=1}^2 x_k\frac{\partial \rho_m}{\partial x_k}m_k = \sum_{m\in S}\rho_m \sum_{m'\in\Delta}\sum_{k=1}^2(\langle m,n_k\rangle - \langle m',n_k\rangle)\rho_{m'}m_k
\]
\[
= \sum_{m\in S}\rho_m\sum_{m'\in\Delta}(m - m')\rho_{m'} = \sum_{m\in S}\rho_m(m - F_{\delta^w}(x))=0.
\]\\
Therefore\\
\[
F_{\delta^w}(x) = \sum_{m\in S}\rho^S_m m = \sum_{m\in S}\frac{|s^{\delta^w} |_m^2}{|s^{\delta^w} |_S^2} m, \mbox{ where } |s^{\delta^w} |_S^2 = \sum_{m\in S}|s^{\delta^w} |_m^2,
\]\\
when $\rho_S = \displaystyle \sum_{m\in S} \rho_m$ achieves maximal. It is easy to derive\\

\begin{tabular}{cl}
\centerline{$\displaystyle F_{\delta^w}(x) = \sum_{m\in S}\rho^S_m m = \sum_{m\not\in S}\rho^{S^c}_m m.$} &  \hspace*{-.4in} $\Box$
\end{tabular}\\
\begin{lm}
\label{de}
There exists a constant $a>0$ (independent of $\delta$) such that for any $x\in P_\Sigma$ the set

\[
S_x = \{m\in \Delta | \rho_m(x) > \delta^a\}
\]

is a simplex in $Z$.
\end{lm}

{\bf Proof:} Take a maximal subset $\tilde{S}_x \subset S_x$ that forms a simplex, which is allowed to contain no integral points in $S_x \setminus \tilde{S}_x$. Clearly, $S_x$ is in the affine span of $\tilde{S}_x$ in $M$. (Without loss of generality, we will assume that $\tilde{S}_x$ forms a top dimensional simplex in $M$. Otherwise, we need to restrict our argument to the affine span of $\tilde{S}_x$ in $M$.) For any $m\in \Delta$, there exists a unique expression

\[
s_m = \delta^{w_m}\prod_{\tilde{m} \in \tilde{S}_x} s_{\tilde{m}}^{l_{\tilde{m}}}.
\]

Correspondingly

\[
\rho_m = \delta^{2w_m}\prod_{\tilde{m} \in \tilde{S}_x} \rho_{\tilde{m}}^{l_{\tilde{m}}}.
\]

For $m \in \Delta$ satisfying $w_m <0$,

\[
\rho_m (x) \geq \delta^{2w_m + a\sum_{\tilde{m} \in \tilde{S}_x} \max(0,l_{\tilde{m}})} > 1
\]

for $a>0$ small. Therefore we may assume $w_m\geq 0$ for all $m\in \Delta$. Since $\{w_m\}_{m\in \Delta}$ is convex and generic, we have $\tilde{S}_x \in Z$. For $m$ not in the simplex spanned by $\tilde{S}_x$, $w_m >0$, we have

\[
\rho_m \leq \delta^{2w_m + a\sum_{\tilde{m} \in \tilde{S}_x} \min(0,l_{\tilde{m}})} \leq \delta^a.
\]

for $a>0$ small. Therefore $S_x = \tilde{S}_x \in Z$.
\hfill$\Box$\\

The following proposition is a direct corollary of lemma \ref{de}.

\begin{prop}
For $S\in Z$ and $x \in P_\Sigma$, assume that $\rho_m (x) > \epsilon$ for all $m\in S$. Then $\rho_m (x) = O(\delta^+)$ for all $m \not\in S$ such that $S\cup\{m\} \not\in Z$.
\hfill$\Box$\\
\end{prop}

{\bf Remark:} In this paper, $O(\delta^+)$ denotes a quantity bounded by $A\delta^a$ for some universal positive constants $A,a$ that only depend on $w$ and $\Delta$. In this paper, the relation between $\epsilon$ and $\delta$ is that we will take $\epsilon$ as small as we want and then take $\delta$ as small as we want depending on $\epsilon$. Geometrically, the metric $\omega_{\delta^w}$ develop necks that have scale $\delta^{a'}$ for some $a' >a$. $\epsilon$ is the gluing scale in section 5 that satisfies $\epsilon \geq \delta^a$. For this section, it is sufficient to take $\epsilon = \delta^a$, which we will assume. In particular, $O(\delta^+) = O(\epsilon)$ in this section.\\

For $S\in Z$, we have 2 $N_{\mathbb{S}}$-invariant \k forms

\[
\omega_{\delta^w}^S = \partial\bar{\partial} h_{\delta^w}^S,\ \ \omega_S = \partial\bar{\partial} h_S,\ \ {\rm where}\ \ h_{\delta^w}^S = \log |s^{\delta^w} |_S^2\ \ h_S = \log |s|_S^2.
\]

The corresponding moment maps are

\[
F_{\delta^w}^S = \sum_{m\in S} \rho^S_m m \mbox{  and  } F_S = \sum_{m\in S} \frac{|s|_m^2}{|s|_S^2} m.
\]

The 2 $N_{\mathbb{S}}$-invariant \k forms and their moment maps coincide if only if $w_m=0$ for $m\in S$.\\

Apply lemma \ref{de}, we have

\begin{prop}
\label{df}
For any $x\in P_\Sigma$, $|\omega_{\delta^w}^{S_x}(x) - \omega_{\delta^w}(x)| = O(\delta^+)$ and $|F_{\delta^w}^{S_x}(x) - F_{\delta^w}(x)| = O(\delta^+)$.
\hfill$\Box$\\
\end{prop}
For each simplex $S\in Z$, let

\[
U^S_\epsilon = \left\{x \in P_\Sigma \left| \rho_S(x) > 1 - |\Delta|\epsilon,\ \rho_m(x) > \epsilon,\ {\rm for}\ m\in S\right.\right\},
\]

where $|\Delta|$ denotes the number of integral points in $\Delta$. The definition clearly implies the following
\begin{prop}
\label{dp}
For any $x\in U^S_\epsilon$, $|\omega_{\delta^w}^{S}(x) - \omega_{\delta^w}(x)| = O(\epsilon)$ and $|F_{\delta^w}^{S}(x) - F_{\delta^w}(x)| = O(\epsilon)$.
\hfill$\Box$\\
\end{prop}
\begin{prop}
\label{cf}
\[
P_\Sigma = \bigcup_{S\in Z} U^S_\epsilon.
\]

Namely, $\{U^S_\epsilon\}_{S\in Z}$ is an open covering of $P_\Sigma$.
\end{prop}
{\bf Proof:} For any $x\in P_\Sigma$, let $S$ contain those $m\in \Delta$ such that $\rho_m(x) > \epsilon$, then $\displaystyle\sum_{m\not\in S}\rho_m(x) \leq |\Delta|\epsilon$. Lemma \ref{de} implies that $S \subset S_x \in Z$ is a simplex. Consequently, $S\in Z$, $x \in U^S_\epsilon$.
\hfill$\Box$\\

Recall $C_{s^{\delta^w}}=(s^{\delta^w})^{-1}(0)$. We have

\begin{prop}
\label{cd}
The image $F_{\delta^w}(C_{s^{\delta^w}})$ is independent of the choice of $w = (w_m)_{m\in \Delta}$ as a representative of an element $[w] \in \tilde{N} = \tilde{N}_0/A$.
\end{prop}

{\bf Proof:} Assume that $\tilde{w} = (\tilde{w}_m)_{m\in \Delta}$ is another representative of $w \in \tilde{N} = \tilde{N}_0/A$. Then there exists $(l,n) \in \mathbb{Z} \oplus N$ such that $\tilde{w}_m = w_m - \langle m,n\rangle + l$. For $x \in N_{\mathbb{C}^*}$, let $\tilde{x} = x+n\log \delta$, then $s_m(\tilde{x}) = \delta^{\langle m,n\rangle} s_m(x)$ and $s^{\delta^{\tilde{w}}}_m (\tilde{x}) = \delta^{\tilde{w}_m} s_m(\tilde{x}) = \delta^l \delta^{w_m} s_m(x) = \delta^l s^{\delta^w}_m(x)$. Hence

\[
s^{\delta^{\tilde{w}}}(\tilde{x}) = \sum_{m\in \Delta} a_m s^{\delta^{\tilde{w}}}_m(\tilde{x}) = \delta^l \sum_{m\in \Delta} a_m s^{\delta^w}_m(x) = \delta^l s^{\delta^w}(x),
\]

and the transformation $x\rightarrow \tilde{x}$ maps $C_{s^{\delta^w}}$ to $C_{s^{\delta^{\tilde{w}}}}$. On the other hand,

\[
|s^{\delta^{\tilde{w}}}_m(\tilde{x})|^2 = \delta^{2l} |s^{\delta^w}_m(x)|^2, \ |s^{\delta^{\tilde{w}}}(\tilde{x})|_\Delta^2 = \delta^{2l} |s^{\delta^w}(x)|_\Delta^2,
\]

\[
F_{\delta^{\tilde{w}}}(\tilde{x}) = \sum_{m\in\Delta} \frac{|s^{\delta^{\tilde{w}}}_m(\tilde{x})|^2}{|s^{\delta^{\tilde{w}}}(\tilde{x})|_\Delta^2}m = \sum_{m\in\Delta} \frac{|s^{\delta^w}_m(x)|^2}{|s^{\delta^w}(x)|_\Delta^2}m = F_{\delta^w}(x).
\]

Therefore $F_{\delta^w}(C_{s^{\delta^w}}) = F_{\delta^{\tilde{w}}}(C_{s^{\delta^{\tilde{w}}}})$.
\hfill $\Box$\\

For each simplex $S\in Z$, let $C_S=s_S^{-1}(0)$, where $s_S = \displaystyle \sum_{m\in S} a_ms_m$, and let $\Gamma_S$ denote the union of all the simplices in the baricenter subdivision of $S$ not containing the vertex of $S$. Then\\
\begin{equation}
\label{cc}
\Gamma_Z = \bigcup_{S\in Z} \Gamma_S
\end{equation}
is a graph in $\Delta$. We have\\
\begin{theorem}
\label{ca}
\[
\lim_{\delta\rightarrow 0}F_{\delta^w}(C_{s^{\delta^w}}) = \bigcup_{S\in Z} F_S(C_S)
\]
is a fattening of $\Gamma_Z$. Consequently, for $\delta\in {\mathbb{R}}_+$ small, $F_{\delta^w}(C_{s^{\delta^w}})$ is a fattening of $\Gamma_Z$.
\end{theorem}

{\bf Proof:} For $x \in P_\Sigma$, according to proposition \ref{cf}, there exists $S \in Z$ such that $x \in U^S_\epsilon$. Since $S$ is a simplex, $w$ can be adjusted by elements in $A$ so that $w_m =0$ for $m\in S$ and $w_m <0$ for $m \not\in S$. According to proposition \ref{cd}, $F_{\delta^w}(C_{s^{\delta^w}})$ is unchanged under such adjustment of $w$. Such adjustment enables us to isolate the discussion to one simplex at a time. For this adjusted weight $w$, $\omega^S_{\delta^w} = \omega_S$ and $F^S_{\delta^w} = F_S$. Proposition \ref{dp} implies that for $x \in U^S_\epsilon$, $F_{\delta^w}(x)$ can be approximated (up to $\epsilon$-terms) by $F_S (x) = F^S_{\delta^w} (x)$.\\

Since $w_m <0$ for $m \not\in S$, we have $|s^{\delta^w} - s_S| = O(\delta^+)$ on $U^S_\epsilon$. $C_{s^{\delta^w}} \cap U^S_\epsilon$ can be approximated (up to $O(\delta^+)$-terms) by $C_S \cap U^S_\epsilon$. Consequently, $F_{\delta^w}(C_{s^{\delta^w}} \cap U^S_\epsilon)$ is an $O(\epsilon)$-approximation of $F_{\delta^w}(C_S \cap U^S_\epsilon)$. Patch such local results together, we get\\
\[
\lim_{\delta\rightarrow 0}F_{\delta^w}(C_{s^{\delta^w}}) = \bigcup_{S\in Z} F_S(C_S).
\]

In fact, $\displaystyle\lim_{\delta\rightarrow 0}C_{s^{\delta^w}} = \bigcup_{S\in Z^{\rm top}} C_S$, where on the righthand side, when $S_1\cap S_2$ is a 1-simplex, the marked points of $C_{S_1}$ and $C_{S_2}$ corresponding to $S_1\cap S_2$ are identified. This limit can be understood in the moduli space $\overline{\cal M}_g$ of stable curves.\\

When $S \in Z$ is a 1-simplex, $F_S(C_S) = \Gamma_S$ is the baricenter of $S$. When $S \in Z$ is a 2-simplex, let $m^0,m^1,m^2$ be the vertices of the simplex $S$. Under the coordinate $x_k = (a_{m^k}s_{m^k})/(a_{m^0}s_{m^0})$ for $k= 1,2$, $C_S = \{x_1 + x_2 + 1\}$ and $F_S(x) = \displaystyle \sum_{k=0}^2 \frac{|x_k|^2}{|x|^2}m^k$, where $x_0=1$ and $|x|^2 = 1 + |x_1|^2 + |x_2|^2$. $F_S(C_S)$ is just the curved triangle in the simplex $S \subset \Delta$ as illustrated in the first picture in figure \ref{fig1}, which is clearly a fattening of the ``Y" shaped graph $\Gamma_S$. Consequently, $\displaystyle \bigcup_{S\in Z} F_S(C_S)$ is a fattening of $\Gamma_Z = \displaystyle \bigcup_{S\in Z} \Gamma_S$.
\hfill$\Box$\\

{\bf Remark:} The result in this theorem is essentially known to Viro in a somewhat different but equivalent form as described in \cite{M}.\\

{\bf Remark:} To achieve the pictures of images of curves in figure 3-5 in the last section, it is necessary to use the moment map introduced in this section. If the moment map of the standard Fubini-Study metric is used, the pictures will look more like hyperbolic metric, more precisely, the holes around center of the polygon will be larger and near the boundary of the polygon will be smaller.\\

{\bf Example:} The Newton polygon\\
\begin{center}
\setlength{\unitlength}{1.1pt}
\begin{picture}(200,100)(-30,10)
\put(150,48){\line(-3,-5){18}}
\put(6,48){\line(3,-5){18}}
\put(132,78){\line(-3,-5){36}}
\put(24,78){\line(3,-5){36}}
\put(42,108){\line(1,0){72}}
\put(114,108){\line(-3,-5){54}}
\put(42,108){\line(3,-5){54}}
\put(24,78){\line(1,0){108}}
\put(78,108){\line(-3,-5){54}}
\put(78,108){\line(3,-5){54}}
\put(6,48){\line(1,0){144}}
\put(42,108){\line(-3,-5){36}}
\put(114,108){\line(3,-5){36}}
\put(24,18){\line(1,0){108}}
\multiput(42,108)(36,0){3}{\circle*{4}}
\multiput(24,78)(36,0){4}{\circle*{4}}
\multiput(6,48)(36,0){5}{\circle*{4}}
\multiput(24,18)(36,0){4}{\circle*{4}}
\end{picture}
\end{center}
\begin{center}
\stepcounter{figure}
Figure \thefigure: the standard simplicial decomposition
\vspace{10pt}
\end{center}
corresponding to an ample line bundle $L$ over $P_\Sigma$ as $\mathbb{P}^2$ with 3 points blown up. Let $E_1,E_2,E_3$ be the 3 exceptional divisors, then\\
\[
L \cong \pi^*({\cal O}(5)) \otimes {\cal O}(-2E_1-E_2-E_3),
\]
where $\pi : P_\Sigma \rightarrow \mathbb{P}^2$ is the natural blow up. Choose a section $s$ of this line bundle near the large complex limit corresponding to the above standard simplicial decomposition of the Newton polygon. Then the curve $C_s = s^{-1}(0)$ cut out by the section $s$ will be mapped to the following under corresponding moment map $F_s$.\\
\begin{center}
\begin{picture}(200,130)(-30,0)
\put(82,128){\makebox(0,0){$E_1$}}
\put(162,23){\makebox(0,0){$E_3$}}
\put(4,23){\makebox(0,0){$E_2$}}
\thicklines
\multiput(-36,-60)(-18,-30){1}
{\put(64,158){\line(4,-3){15.2}}
\put(79.2,146.6){\line(0,-1){17}}
\put(46,128){\line(3,-1){18}}
\put(64,122){\line(2,1){15.4}}}

\multiput(-54,-90)(-18,-30){1}
{\put(64,158){\line(4,-3){15.2}}
\put(79.2,146.6){\line(0,-1){16.5}}
\put(64,119){\line(4,3){15.2}}}

\multiput(200,-60)(18,-30){1}
{\put(-64,158){\line(-4,-3){15.2}}
\put(-79.2,146.6){\line(0,-1){17}}
\put(-46,128){\line(-3,-1){18}}
\put(-64,122){\line(-2,1){15.4}}}

\multiput(218,-90)(18,-30){1}
{\put(-64,158){\line(-4,-3){15.2}}
\put(-79.2,146.6){\line(0,-1){16.5}}
\put(-64,119){\line(-4,3){15.2}}}

\multiput(-18,-90)(36,0){2}
{\put(64,98){\line(1,6){3.3}}
\put(100,98){\line(-1,6){3.3}}
\put(82,125){\line(2,-1){15}}
\put(82,125){\line(-2,-1){15}}}

\multiput(0,217)(36,0){1}
{\put(64,-98){\line(1,-6){3.3}}
\put(100,-98){\line(-1,-6){3.3}}
\put(82,-125){\line(2,1){15}}
\put(82,-125){\line(-2,1){15}}}

\multiput(-54,-90)(36,0){1}
{\put(100,98){\line(-1,6){3.3}}
\put(82,125){\line(2,-1){15}}
\put(82,125){\line(-3,-1){18}}}

\multiput(-36,217)(36,0){1}
{\put(100,-98){\line(-1,-6){3.3}}
\put(82,-125){\line(2,1){15}}
\put(82,-125){\line(-3,1){18}}}

\multiput(218,-90)(36,0){1}
{\put(-100,98){\line(1,6){3.3}}
\put(-82,125){\line(-2,-1){15}}
\put(-82,125){\line(3,-1){18}}}

\multiput(200,217)(36,0){1}
{\put(-100,-98){\line(1,-6){3.3}}
\put(-82,-125){\line(-2,1){15}}
\put(-82,-125){\line(3,1){18}}}

\multiput(100,78.5)(18,-30){2}{\circle{30}}
\multiput(64,78.5)(18,-30){2}{\circle{30}}
\multiput(46,48.5)(18,-30){1}{\circle{30}}

\thinlines
\multiput(64,98)(36,0){2}{\line(0,1){21}}
\multiput(28,98)(36,0){3}{\line(2,-1){18}}
\multiput(64,98)(36,0){3}{\line(-2,-1){18}}
\multiput(46,68)(36,0){3}{\line(0,1){21}}
\multiput(10,68)(36,0){4}{\line(2,-1){18}}
\multiput(46,68)(36,0){4}{\line(-2,-1){18}}
\multiput(28,38)(36,0){4}{\line(0,1){21}}
\multiput(28,38)(36,0){4}{\line(2,-1){18}}
\multiput(28,38)(36,0){4}{\line(-2,-1){18}}
\multiput(46,8)(36,0){3}{\line(0,1){21}}
\put(40.6,119){\line(1,0){82.8}}
\put(22.25,8){\line(1,0){119.5}}
\put(22.25,8){\line(-3,5){24}}
\put(141.75,8){\line(3,5){24}}
\put(40.6,119){\line(-3,-5){42.5}}
\put(123.4,119){\line(3,-5){42.5}}
\end{picture}
\end{center}
\begin{center}
\refstepcounter{figure} \label{fig3}
Figure \thefigure: $F_s(C_s)$
\vspace{10pt}
\end{center}

\subsection{Secondary fan}
Theorem \ref{ca} can be better understood in the context of the secondary fan. To begin with, we consider the space ${\cal M}_\Delta$ of curves $C_s$ modulo the equivalent relations of toric actions. With a little abuse of notation, we will call ${\cal M}_\Delta$ the toric moduli space of curves $C_s$ with the Newton polygon $\Delta$. Let $\tilde{M}_0 \cong \mathbb{Z}^{\Delta}$ be the dual lattice of $\tilde{N}_0 = \{ w=(w_m)_{m\in \Delta}\in \mathbb{Z}^{\Delta} \} \cong \mathbb{Z}^{\Delta}$.\\

Recall that $\tilde{N} = \tilde{N}_0/A$. The dual lattice $\tilde{M}=A^\perp$. We have the natural identification

\[
{\cal M}_\Delta \cong \tilde{N}_{\mathbb{C}^*} ={\rm Spec}(\mathbb{C}[\tilde{M}]) \cong (\mathbb{C}^*)^{\Delta}/N^+_{\mathbb{C}^*}.
\]

To make sense of the large complex limit, we need the compactification $\overline{{\cal M}_\Delta}$ of ${\cal M}_\Delta$ determined by the so-called secondary fan.\\

For general $[w] \in \tilde{N}$, $w=(w_m)_{m\in \Delta}$ is not convex on $\Delta$. Let $\hat{w} =(\hat{w}_m)_{m\in \Delta}$ be the convex hull of $w$. When $w$ is generic, $\hat{w}$ determines a simplicial decomposition $Z_w$ of $\Delta$. (It is easy to observe that $Z_w$ is independent of the choice of representative $w$ in the equivalent class $[w]$.) Let $\hat{S}$ be the set of $S\subset \Delta$ that forms an $r$-dimensional simplex. Then $Z_w$ can be regarded as a subset of $\hat{S}$. Let $\hat{Z}$ denote the set of all $Z_w$ for $[w] \in \tilde{N}$. For $Z \in \hat{Z}$, let $\tau_Z \subset \tilde{N}$ be the closure of the set of all $[w] \in \tilde{N}$ such that $Z_w = Z$. Each $\tau_Z$ is a convex integral top dimension cone in $\tilde{N}$. The union of all $\tau_Z$ is exactly $\tilde{N}$. Let $\hat{\Sigma}$ be the fan whose cones are subcones of the top dimensional cones $\{\tau_Z\}_{Z\in \hat{Z}}$. $\hat{\Sigma}$ is a complete fan.\\

Let $\tilde{Z}$ be the set of simplicial decompositions $Z_w\subset \tilde{S}$ of $\Delta$ that is determined by some strictly convex function $w=(w_m)_{m\in \Delta}$ on $\Delta$. For $Z \in \tilde{Z}$, let $\tau_Z \subset \tilde{N}$ be the set of $[w] \in \tilde{N}$, where $w=(w_m)_{m\in \Delta}$ is a piecewise linear convex function on $\Delta$ with respect to the simplicial decomposition $Z$. Each $\tau_Z$ is a integral top dimension cone in $\tilde{N}$. The union of all $\tau_Z$\\
\[
\tau = \bigcup_{Z\in \tilde{Z}} \tau_Z
\]\\
is exactly the convex cone of all $[w] \in \tilde{N}$, where $w=(w_m)_{m\in \Delta}$ is a piecewise linear convex function on $\Delta$. Let $\tilde{\Sigma}$ be the fan whose cones are subcones of the top dimensional cones $\{\tau_Z\}_{Z\in \tilde{Z}}$. $\tilde{\Sigma}$ is a subfan of the complete fan $\hat{\Sigma}$.\\

The fan $\hat{\Sigma}$ is the so-called {\bf secondary fan}. (For more detail about the secondary fan, please refer to the book \cite{G}. \cite{AGM} contains some application of secondary fan to mirror symmetry.) $\hat{\Sigma}$ naturally determines the compactification $\overline{{\cal M}_\Delta} = P_{\hat{\Sigma}}$. We will call $\tilde{\Sigma}$ the {\bf partial secondary fan}. $\tilde{\Sigma}$ determines the partial compactification $\widetilde{{\cal M}_\Delta} = P_{\tilde{\Sigma}}$. For each $Z\in \hat{Z}$, the top dimensional cone $\tau_Z$ determines a single fixed point $s_\infty^Z\in \overline{{\cal M}_\Delta}\backslash {\cal M}_\Delta$ of the $\tilde{N}_{\mathbb{C}^*}$ action. We will call such $s_\infty^Z$ a {\bf large complex limit} point. The set of different large complex limit points is parameterized by the set of simplicial decomposition $\hat{Z}$. Each large complex limit point $s_\infty^Z$ possesses a cell neighborhood $\overline{\tau_Z^{\mathbb{C}}} \subset \overline{{\cal M}_\Delta}$, where $\tau_Z^{\mathbb{C}} = \tau_Z \otimes _{\mathbb{Z}_{\geq 0}} \mathbb{C}_+ \subset \tilde{N}_{\mathbb{C}^*}$, $\mathbb{Z}_{\geq 0}$ acts trivially on $\mathbb{C}_+ = \{z\in \mathbb{C}^*: |z| \geq 1\}$. We have the following natural cell decomposition of $\overline{{\cal M}_\Delta}$\\
\[
\overline{{\cal M}_\Delta} =\bigcup_{Z\in \tilde{Z}} \overline{\tau_Z^{\mathbb{C}}}
\]\\
Given a simplicial decomposition $Z\in \hat{Z}$ of $\Delta$, let $\tau^0_Z$ denote the interior of $\tau_Z$. Any $[w]\in \tau^0_Z$ can be represented by a strongly convex piecewise linear function $w=(w_m)_{m\in \Delta}$ on $\Delta$ with respect to $Z$. It is easy to see that when $\delta$ approaches 0, $C_{s^{\delta^w}}$ will approach the large complex limit point $s_\infty^Z$ in $\overline{{\cal M}_\Delta}$. In such situation, we will say that $C_{s^{\delta^w}}$ or $s_\delta$ is {\bf near the large complex limit} point (determined by $Z$), when $\delta$ is small.\\

Theorem \ref{ca} applies to each of such large complex limit point $s_\infty^Z$ in $\widetilde{{\cal M}_\Delta}$ for $Z \in \tilde{Z}$, and can be rephrased as: when the string diagrams $C_{s^{\delta^w}}$ approach the large complex limit point $s_\infty^Z$ in $\widetilde{{\cal M}_\Delta}$ as $\delta \rightarrow 0$, the amoebas $F_{\delta^w}(C_{s^{\delta^w}})$ of the string diagrams $C_{s^{\delta^w}}$ converge to the Feynman diagram $\Gamma_Z$.\\

Theorem \ref{ca} can be generalized to the full compactification $\overline{{\cal M}_\Delta}$ without additional difficulty.\\
\begin{theorem}
\label{cg}
For $Z \in \hat{Z}$, when the string diagrams $C_{s^{\delta^w}}$ approach the large complex limit point $s_\infty^Z$ in $\overline{{\cal M}_\Delta}$ as $\delta \rightarrow 0$, the amoebas $F_{\delta^w}(C_{s^{\delta^w}})$ of the string diagrams $C_{s^{\delta^w}}$ converge to the Feynman diagram $\Gamma_Z$.
\end{theorem}
{\bf Proof:} It is straightforward to generalize lemma \ref{de}, propositions \ref{df}, \ref{dp}, \ref{cf}, \ref{cd} and in particular, theorem \ref{ca} to the case when $Z \in \hat{Z}$. The arguments are literally the same with the understanding that $S$ considered as a subset in $M$ contains only the integral vertex points of the simplex $S$, not any other integral points in the simplex $S$.
\hfill$\Box$\\

{\bf Remark:} Theorem \ref{ca} is used in \cite{lag3} to construct Lagrangian torus fibration for quintic Calabi-Yau manifolds near large complex limit in the partial secondary fan compactification. Theorem \ref{cg} can be used to construct similar Lagrangian torus fibration for quintic Calabi-Yau manifolds near large complex limit that is not necessarily in the partial secondary fan compactification. According to \cite{AGM}, a large complex limit in the partial secondary fan compactification, under the mirror symmetry, corresponds to large radius limit of a \k cone of the mirror Calabi-Yau manifold, while a large complex limit not in the partial secondary fan compactification, under the mirror symmetry, may correspond to large radius limit of some other physical model like Landau-Ginzberg model etc.\\

\se{The pair of pants and the 3-valent vertex of a graph}
In Feynman diagram, a 3-valent vertex represents the most basic particle interaction. In string theory, the corresponding string diagram is the pair of pants, which can be represented by a general line in ${\mathbb{CP}^2}$ with the 3 punctured points being the intersection points of this line with the 3 coordinate lines. In this section, we will describe an analogue of this picture in our situation. More precisely, the standard moment map maps a general line to a fattening of the 3-valent vertex neighborhood of a graph. In this section, we will explicitly perturb the moment map, so that the perturbed moment map will map the general line to the 3-valent vertex neighborhood, i.e., a ``Y" shaped graph.\\

\subsection{The piecewise smooth case}
Consider ${\mathbb{CP}^2}$ with the Fubini-Study metric and the curve $C_0: z_0 + z_1 + z_2 =0$ in ${\mathbb{CP}^2}$. We have the torus fibration $F: {\mathbb{CP}^2}\rightarrow \mathbb{R}^+\mathbb{P}^2$ defined as
\[
F( [z_1,z_2,z_3]) = [|z_1|,|z_2|,|z_3|].
\]
Under the inhomogenuous coordinate $x_i= z_i/z_0$, locally we have
\[
F: \mathbb{C}^2\rightarrow (\mathbb{R}^+)^2, F(x_1,x_2)=(r_1,r_2),
\]
where $x_k=r_ke^{i\theta_k}$. The image of $C_0: x_1 + x_2 +1=0$ under $F$ is\\
\[
\tilde{\Gamma} =\{(r_1,r_2)|r_1 + r_2 \geq 1, r_1 \leq r_2 + 1, r_2 \leq r_1 + 1 \}.
\]
$C_0$ is a symplectic submanifold. We want to deform $C_0$ symplectically to $C_1$ whose image under $F$ is expected to be\\
\[
\Gamma = \{(r_1,r_2)| 0\leq r_2\leq r_1=1 \ {\rm or} \ 0\leq r_1\leq r_2 =1 \ {\rm or} \ r_1=r_2\geq 1\}.
\]
A moment of thought suggests taking $C_t = {\cal F}_t (C_0)$, where\\
\[
{\cal F}_t (x_1, x_2) = \left(\textstyle \left(\frac{\max(1,r_2)}{\max(r_1,r_2)}\right)^tx_1, \left(\frac{\max(1,r_1)}{\max(r_1,r_2)}\right)^tx_2\right).
\]\\
The \k form of the Fubini-Study metric can be written as\\
\[
\omega_{\rm FS} = \frac{ dx_1\wedge d\bar{x}_1 + dx_2\wedge d\bar{x}_2 + (x_2 dx_1 - x_1 dx_2)\wedge (\bar{x}_2d\bar{x}_1-\bar{x}_1d\bar{x}_2)}{(1+|x|^2)^2}.
\]

\begin{lm}
\label{da}
$\omega_{\rm FS}$ restricts to a symplectic form on $C_t \setminus {\rm Sing}(C_t)$, where ${\rm Sing}(C_t) := \{x \in C_t: (r_1-1)(r_2-1)(r_1-r_2) = 0\}$. More precisely, there exists $c>1$ such that $\frac{1}{c} \omega_{\rm FS} \leq {\cal F}_t^* \omega_{\rm FS} \leq c \omega_{\rm FS}$ on $C_0$ for all $t \in [0,1]$.
\end{lm}
{\bf Proof:} Due to the symmetries of permuting $[z_0,z_1,z_2]$, to verify that $C_t$ is symplectic, we only need to verify for one region out of six. Consider $1\geq |x_2|\geq |x_1|$, where\\
\[
C_t = \left\{\textstyle \left( \left(\frac{1}{r_2}\right)^tx_1,\left(\frac{1}{r_2}\right)^tx_2\right): x_1 + x_2 +1=0\right\}.
\]
$x_1 + x_2 +1=0$ implies that\\
\[
dx_1 = - dx_2.
\]
Recall that\\
\[
\textstyle \frac{dr_k}{r_k}= {\rm Re}\left(\frac{dx_k}{x_k}\right),\ d\theta_k = {\rm Im}\left(\frac{dx_k}{x_k}\right).
\]
Consequently\\
\[
\textstyle \frac{dr_1}{r_1}= {\rm Re}\left(\frac{dx_1}{x_1}\right)= -{\rm Re}\left(\left(\frac{x_2}{x_1}\right)\frac{dx_2}{x_2}\right),
\]
\[
\textstyle \frac{dr_2}{r_2}= {\rm Re}\left(\frac{dx_2}{x_2}\right)= -{\rm Re}\left(\left(\frac{x_1}{x_2}\right)\frac{dx_1}{x_1}\right).
\]\\
We have\\
\[
\textstyle d\left(\left(\frac{1}{r_2}\right)^tx_1\right) = \left(\frac{1}{r_2}\right)^t\left(dx_1 - tx_1\frac{dr_2}{r_2}\right)
\]
\[
\textstyle d\left(\left(\frac{1}{r_2}\right)^tx_1\right)\wedge d\left(\left(\frac{1}{r_2}\right)^t\bar{x}_1\right) = \left(\frac{1}{r_2}\right)^{2t}\left(dx_1\wedge d\bar{x}_1 +t(x_1d\bar{x}_1 - \bar{x}_1dx_1)\wedge \frac{dr_2}{r_2}\right)
\]
\[
\textstyle =\left(\frac{1}{r_2}\right)^{2t} \left(1 + t{\rm Re}\left(\frac{x_1}{x_2}\right)\right)dx_1\wedge d\bar{x}_1,
\]\\
\[
\textstyle d\left(\left(\frac{1}{r_2}\right)^tx_2\right) = \left(\frac{1}{r_2}\right)^t\left(dx_2 - tx_2\frac{dr_2}{r_2}\right)
\]
\[
\textstyle d\left(\left(\frac{1}{r_2}\right)^tx_2\right)\wedge d\left(\left(\frac{1}{r_2}\right)^t\bar{x}_2\right)
\]
\[
\textstyle = \left(\frac{1}{r_2}\right)^{2t}\left(dx_2\wedge d\bar{x}_2 +t(x_2d\bar{x}_2 - \bar{x}_2dx_2)\wedge \frac{dr_2}{r_2}\right) =\left(\frac{1}{r_2}\right)^{2t} (1 - t)dx_2\wedge d\bar{x}_2,
\]\\
\[
\textstyle \left( \left(\frac{1}{r_2}\right)^tx_2\right)d\left( \left(\frac{1}{r_2}\right)^tx_1\right) - \left( \left(\frac{1}{r_2}\right)^tx_1\right)d\left( \left(\frac{1}{r_2}\right)^tx_2\right)
\]
\[
\textstyle = \left(\frac{1}{r_2}\right)^{2t}(x_2dx_1 - x_1dx_2) = \left(\frac{1}{r_2}\right)^{2t}x_2\left(1+ \left(\frac{x_1}{x_2}\right)\right)dx_1 = -\left(\frac{1}{r_2}\right)^{2t} dx_1.
\]\\
By restriction to $C_t$ and use the fact that $1 + {\rm Re}\left(\frac{x_1}{x_2}\right) \geq \frac{1}{2}$ on $C_0$, we get\\
\[
\frac{({\cal F}_t^* \omega_{\rm FS})|_{C_0}}{dx_1\wedge d\bar{x}_1} = \frac{(1-t)\left(\frac{1}{r_2}\right)^{2t} + \left(\frac{1}{r_2}\right)^{2t} \left(1 + t{\rm Re}\left(\frac{x_1}{x_2}\right)\right) + \left(\frac{1}{r_2}\right)^{4t}}{\left(1 + \left(\frac{1}{r_2}\right)^{2t} (r_1^2 + r_2^2) \right)^2} \geq \frac{1}{6}.
\]\\
\[
\frac{({\cal F}_t^* \omega_{\rm FS})|_{C_0}}{\omega_{\rm FS}|_{C_0}} = \frac{(1-t)\left(\frac{1}{r_2}\right)^{2t} + \left(\frac{1}{r_2}\right)^{2t} \left(1 + t{\rm Re}\left(\frac{x_1}{x_2}\right)\right) + \left(\frac{1}{r_2}\right)^{4t}}{3\left(\frac{1}{1 + r_1^2 + r_2^2} + \left(\frac{1}{r_2}\right)^{2t} \frac{r_1^2 + r_2^2}{1 + r_1^2 + r_2^2} \right)^2} \geq \frac{1}{2}.
\]\\
These computations show that $C_t$ is symplectic in the region $r_1 < r_2 < 1$. By symmetry, we can see that $C_t$ is symplectic in the other five regions.
\hfill$\Box$\\
\begin{prop}
\label{dm}
${\cal F}_t^* \omega_{\rm FS}$ is a piecewise smooth continuous symplectic form on $C_0$ for any $t \in [0,1]$.
\end{prop}
{\bf Proof:} In light of lemma \ref{da}, only continuity need comment. This is an easy consequence of the invariance of ${\cal F}_t^* \omega_{\rm FS}$ under the symmetries of mutating the coordinate $[z_0,z_1,z_2]$.
\hfill$\Box$\\
\begin{theorem}
\label{dc}
There exists a family of piecewise smooth Lipschitz Hamiltonian diffeomorphism $H_t: {\mathbb{CP}^2}\rightarrow{\mathbb{CP}^2}$ such that $H_t$ is smooth away from ${\rm Sing}(C_0)$, $H_t(C_0)=C_t$, $H_t({\rm Sing}(C_0))= {\rm Sing}(C_t)$ and $H_t$ is identity away from an arbitrary small neighborhood of $C_{[0,1]} := \displaystyle \bigcup_{t \in [0,1]}C_t$. In particular $H_t$ leaves $\partial \mathbb{CP}^2$ (the union of the three coordinate ${\mathbb{CP}^1}$'s) invariant. The perturbed moment map (Lagrangian fibration) $\hat{F} = F\circ H_1$ satisfies $\hat{F}(C_0) = \Gamma$ (the ``Y" shaped graph with a 3-valent vertex $v_0$).\\
\end{theorem}
{\bf Proof:} Lemma \ref{da} implies that $C_t$'s are piecewise smooth symplectic submanifolds in ${\mathbb{CP}^2}$. Each $C_t$ is a union of 6 pieces of smooth symplectic submanifolds with boundaries and corners. The 6 pieces have equal area (equal to one-sixth of the total area of $C_t$), which is independent of $t$. $C_0$ is symplectic isotopic to $C_1$ via the family $\{C_t\}$. By extension theorem (corollary 6.3) in \cite{lag2}, we may construct a piecewise smooth Lipschitz Hamiltonian diffeomorphism $H_t: {\mathbb{CP}^2}\rightarrow{\mathbb{CP}^2}$ such that $H_t(C_0)=C_t$. Corollary 6.3 in \cite{lag2} can further ensure that $H_t$ leaves $\partial \mathbb{CP}^2$ invariant as desired.\\

More precisely, the proof of corollary 6.3 in \cite{lag2} is separated into 2 steps. In the first step, one modify the symplectic isotopy (see section 6 of \cite{lag2} for definition) ${\cal F}_t: C_0 \rightarrow C_t$ into a symplectic flow while keeping the restriction of ${\cal F}_t$ to the boundaries of the 6 pieces unchanged. (one in fact first modify ${\cal F}_t$ in one of the 6 pieces, then extend the modification symmetrically to the other pieces.) In particular, $C_t \cap \partial \mathbb{CP}^2$ is fixed by the symplectic flow. In the second step, theorem 6.9 in \cite{lag2} is applied to extend the symplectic flow to $\mathbb{CP}^2$ while keeping $\partial \mathbb{CP}^2$ fixed. The construction in effect ensures that $H_t|_{{\rm Sing}(C_0)} = {\cal F}_t|_{{\rm Sing}(C_0)}$ and $H_t$ is smooth away from ${\rm Sing}(C_0)$.
\hfill$\Box$\\

Similar construction can be carried out for degree $d$ Fermat type curves. (The case of $d=5$ is carried out in \cite{lag2}.)\\
\begin{center}
\begin{picture}(200,70)(-50,120)
\put(-100,0){\thicklines
\put(64,158){\line(3,-1){12}}
\put(76,154){\line(1,0){12}}
\put(100,158){\line(-3,-1){12}}

\put(18,30)
{\put(46,128){\line(4,-3){10}}
\put(56,120.5){\line(3,-5){6}}
\put(64,98){\line(-1,6){2.1}}}

\put(146,30)
{\put(-46,128){\line(-4,-3){10}}
\put(-56,120.5){\line(-3,-5){6}}
\put(-64,98){\line(1,6){2.1}}}

\thinlines
\multiput(82,149)(36,0){1}{\line(2,1){18}}
\multiput(82,149)(36,0){1}{\line(-2,1){18}}
\multiput(82,128)(36,0){1}{\line(0,1){21}}
\put(46,128){\line(1,0){72}}
\put(46,128){\line(3,5){36}}
\put(118,128){\line(-3,5){36}}}

\thicklines
\put(30,158){\vector(1,0){40}}

\put(180,0){\thicklines
\multiput(-68,149)(36,0){1}{\line(2,1){18}}
\multiput(-68,149)(36,0){1}{\line(-2,1){18}}
\multiput(-68,128)(36,0){1}{\line(0,1){21}}

\thinlines
\put(-104,128){\line(1,0){72}}
\put(-104,128){\line(3,5){36}}
\put(-32,128){\line(-3,5){36}}}
\end{picture}
\end{center}
\begin{center}
\stepcounter{figure}
Figure \thefigure: $F(C_p)$ of $p=z_1^d + z_2^d + z_3^d$ perturbed to $\hat{F}(C_p) = \Gamma$
\end{center}
Let $\mu_r := C_0 \cap \hat{F}^{-1}(r)$ for $r\in \Gamma$. When $d=1$, for $r$ being one of the three boundary points of $\Gamma$, $\mu_r$ is a point. For $r$ in smooth part of $\Gamma$, $\mu_r$ is a circle. For $r$ being the unique singular point of $\Gamma$, which in quantum mechanics usually indicate the particle interaction point, $\mu_r$ is of ``$\Theta$'' shape. This picture indicates the simplest string interaction.\\\\
When $d=5$, for $r$ being one of the three boundary points of $\Gamma$, $\mu_r$ is 5 points. For $r$ in smooth part of $\Gamma$, $\mu_r$ is 5 circles. For $r$ being the unique singular point of $\Gamma$, which in quantum mechanics usually indicate the particle interaction point, $\mu_r$ is a graph in 2-torus as indicated in the following picture, which is much more complicated than $d=1$ case. This picture indicates sort of degenerate multi-particle string interaction with multiplicity.\\\\
\begin{center}
\leavevmode
\hbox{%
\epsfxsize=4in
\epsffile{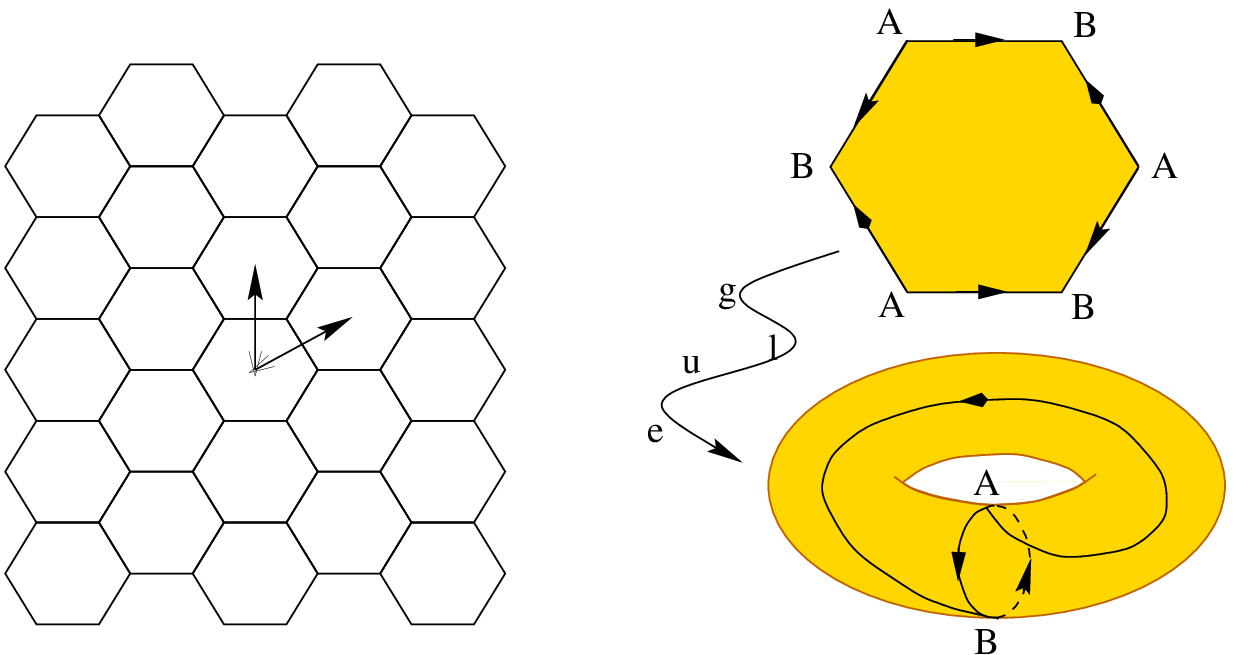}}
\end{center}
\begin{center}
\refstepcounter{figure} \label{fig9}
Figure \thefigure: $\hat{F}^{-1}({\rm Sing}(\Gamma))$ for $d=5$ and $d=1$
\vspace{10pt}
\end{center}

\subsection{The smooth case}
Notice that $C_t$ in section 4.1 is not smooth on ${\rm Sing}(C_t)$. In this section, we will make $C_t$ smooth. The trade-off is that $F(C_1) = \Gamma$ except in a small neighborhood of the vertex of the graph $\Gamma$, where $F(C_1)$ is a fattening of $\Gamma$. To modify the definition of $C_t$ to make it smooth, consider real function $h(a)\geq 0$ such that $h(a)+h(-a)=1$ for all $a$ and $h(a)=0$ for $a\leq -\epsilon$. Then consequently, $h(a)=1$ for $a\geq \epsilon$ and $h(a)\leq 1$.\\

We may modify the definition of $C_t$ to consider $C_t = \tilde{\cal F}_t(C_0)$, where

\[
\tilde{\cal F}_t (x_1,x_2) = \left(\textstyle \left(\frac{\eta_1}{\eta_0}\right)^tx_1, \left(\frac{\eta_2}{\eta_0}\right)^tx_2\right),
\]
\[
\eta_2 = r_1^{h(\log r_1)},\ \eta_1 = r_2^{h(\log r_2)},\ \eta_0 = r_1\textstyle \left(\frac{r_2}{r_1}\right)^{h(\log (r_2/r_1))}.
\]

$C_t$ is now smooth and is only modified in a $\epsilon$-neighborhood of ${\rm Sing}(C_t)$.\\

Assume $\lambda(a) = h(a) + h'(a)a$, $\lambda_0 = \lambda(\log r_2 - \log r_1)$, $\lambda_1 = \lambda(\log r_1)$, $\lambda_2 = \lambda(\log r_2)$. Then\\
\[
\textstyle \frac{d\eta_2}{\eta_2} = \lambda_1\frac{dr_1}{r_1},\ \frac{d\eta_1}{\eta_1} = \lambda_2\frac{dr_2}{r_2},\ \frac{d\eta_0}{\eta_0} = \frac{dr_1}{r_1} + \lambda_0\left(\frac{dr_2}{r_2} - \frac{dr_1}{r_1}\right).
\]\\
\[
\textstyle d\left(\left(\frac{\eta_1}{\eta_0}\right)^tx_1\right) = \left(\frac{\eta_1}{\eta_0}\right)^t\left(dx_1 + tx_1\left(\frac{d\eta_1}{\eta_1} - \frac{d\eta_0}{\eta_0}\right)\right).
\]
\[
\textstyle \frac{d\eta_1}{\eta_1} - \frac{d\eta_0}{\eta_0} = -(1 - \lambda_0) \frac{dr_1}{r_1} - (\lambda_0 - \lambda_2)\frac{dr_2}{r_2}.
\]
\begin{eqnarray*}
&&\textstyle d\left(\left(\frac{\eta_1}{\eta_0}\right)^tx_1\right)\wedge d\left(\left(\frac{\eta_1}{\eta_0}\right)^t\bar{x}_1\right)\\
&=& \textstyle \left(\frac{\eta_1}{\eta_0}\right)^{2t}\left(dx_1\wedge d\bar{x}_1 -t(x_1d\bar{x}_1 - \bar{x}_1dx_1)\wedge \left(\frac{d\eta_1}{\eta_1} - \frac{d\eta_0}{\eta_0}\right)\right)\\
&=&\textstyle \left(\frac{\eta_1}{\eta_0}\right)^{2t} \left(1 - (1-\lambda_0)t + (\lambda_0 - \lambda_2)t{\rm Re}\left(\frac{x_1}{x_2}\right)\right)dx_1\wedge d\bar{x}_1.
\end{eqnarray*}

\[
\textstyle d\left(\left(\frac{\eta_2}{\eta_0}\right)^tx_2\right) = \left(\frac{\eta_2}{\eta_0}\right)^t\left(dx_2 + tx_2\left(\frac{d\eta_2}{\eta_2} - \frac{d\eta_0}{\eta_0}\right)\right).
\]
\[
\textstyle \frac{d\eta_2}{\eta_2} - \frac{d\eta_0}{\eta_0} = -(1 - \lambda_0 - \lambda_1) \frac{dr_1}{r_1} - \lambda_0 \frac{dr_2}{r_2}.
\]
\begin{eqnarray*}
&&\textstyle d\left(\left(\frac{\eta_2}{\eta_0}\right)^tx_2\right)\wedge d\left(\left(\frac{\eta_2}{\eta_0}\right)^t\bar{x}_2\right)\\
&=& \textstyle \left(\frac{\eta_2}{\eta_0}\right)^{2t}\left(dx_2\wedge d\bar{x}_2 +t(x_2d\bar{x}_2 - \bar{x}_2dx_2)\wedge \left(\frac{d\eta_2}{\eta_2} - \frac{d\eta_0}{\eta_0}\right)\right)\\
&=&\textstyle \left(\frac{\eta_2}{\eta_0}\right)^{2t} \left(1+ (1-\lambda_0 - \lambda_1)t{\rm Re}\left(\frac{x_2}{x_1}\right) - \lambda_0 t\right)dx_2\wedge d\bar{x}_2.
\end{eqnarray*}

\[
\textstyle \alpha = \left( \left(\frac{\eta_2}{\eta_0}\right)^tx_2\right)d\left( \left(\frac{\eta_1}{\eta_0}\right)^tx_1\right) - \left( \left(\frac{\eta_1}{\eta_0}\right)^tx_1\right)d\left( \left(\frac{\eta_2}{\eta_0}\right)^tx_2\right)
\]
\[
\textstyle = \left(\frac{\eta_2\eta_1}{\eta_0^2}\right)^{t}\left(x_2dx_1 - x_1dx_2 + tx_1x_2\left(\frac{d\eta_1}{\eta_1} - \frac{d\eta_2}{\eta_2}\right)\right)
\]
\[
\textstyle = \left(\frac{\eta_2\eta_1}{\eta_0^2}\right)^{t} \left(-dx_1 + tx_1x_2\left(\lambda_2\frac{dr_2}{r_2} - \lambda_1\frac{dr_1}{r_1}\right)\right).
\]\\
\[
\textstyle \alpha \wedge \bar{\alpha} = \left(\frac{\eta_2\eta_1}{\eta_0^2}\right)^{2t}\left(dx_1d\bar{x}_1 - t\left(\bar{x}_1\bar{x}_2dx_1 - x_1x_2d\bar{x}_1\right)\left(\lambda_2\frac{dr_2}{r_2} - \lambda_1\frac{dr_1}{r_1}\right)\right)
\]
\[
\textstyle = \left(\frac{\eta_2\eta_1}{\eta_0^2}\right)^{2t}\left(1 + t\left(\lambda_2{\rm Re}(x_1) + \lambda_1{\rm Re}(x_2)\right)\right)dx_1d\bar{x}_1.
\]
By restriction to $C_t$ we get\\
\[
\frac{(\tilde{\cal F}_t^* \omega_{\rm FS})|_{C_0}}{dx_1\wedge d\bar{x}_1}= \textstyle \left[\left(\frac{\eta_2}{\eta_0}\right)^{2t}\left(1+ (1-\lambda_0 - \lambda_1)t{\rm Re}\left(\frac{x_2}{x_1}\right) - \lambda_0 t\right)\right.
\]
\[
\textstyle + \left(\frac{\eta_1}{\eta_0}\right)^{2t}\left(1 - (1-\lambda_0)t + (\lambda_0 - \lambda_2)t{\rm Re}\left(\frac{x_1}{x_2}\right)\right)
\]
\[
\textstyle + \left.\left.\left(\frac{\eta_2\eta_1}{\eta_0^2}\right)^{2t}\left(1 + t\left(\lambda_2{\rm Re}(x_1) + \lambda_1{\rm Re}(x_2)\right)\right)\right]\right/\left(1+\left(\frac{\eta_2}{\eta_0}\right)^{2t}r_2^2 + \left(\frac{\eta_1}{\eta_0}\right)^{2t}r_1^2\right)^2
\]
\[
= \frac{\tilde{\omega}_t}{dx_1\wedge d\bar{x}_1} + tR_t, \mbox{ where } R_t = (1- \lambda_0) R_{t,0} + \lambda_1 R_{t,1} + \lambda_2 R_{t,2},
\]
\[
R_{t,0} = \frac{\left(\frac{\eta_2}{\eta_0}\right)^{2t}\left(1+ {\rm Re}\left(\frac{x_2}{x_1}\right)\right) - \left(\frac{\eta_1}{\eta_0}\right)^{2t}\left(1 + {\rm Re}\left(\frac{x_1}{x_2}\right)\right)}{\left(1+ \left(\frac{\eta_2}{\eta_0} \right)^{2t} r_2^2 + \left(\frac{\eta_1}{\eta_0}\right)^{2t}r_1^2\right)^2}
\]
\[
R_{t,1} = \frac{\left(\frac{\eta_2}{\eta_0}\right)^{2t}\left(\eta_1^{2t} {\rm Re}(x_2) - {\rm Re}\left(\frac{x_2}{x_1}\right)\right)}{\left(1+ \left(\frac{\eta_2}{\eta_0} \right)^{2t} r_2^2 + \left(\frac{\eta_1}{\eta_0}\right)^{2t}r_1^2\right)^2},\ R_{t,2} = \frac{\left(\frac{\eta_1}{\eta_0}\right)^{2t}\left(\eta_2^{2t} {\rm Re}(x_1) - {\rm Re} \left(\frac{x_1}{x_2}\right)\right)}{\left(1+ \left(\frac{\eta_2}{\eta_0} \right)^{2t} r_2^2 + \left(\frac{\eta_1}{\eta_0}\right)^{2t}r_1^2\right)^2},
\]
\[
\tilde{\omega}_t = \frac{\left(\frac{\eta_2}{\eta_0}\right)^{2t}\left(1+ t{\rm Re}\left(\frac{x_1}{x_2}\right)\right) + \left(\frac{\eta_1}{\eta_0}\right)^{2t}(1 - t) + \left(\frac{\eta_2\eta_1}{\eta_0^2}\right)^{2t}}{\left(1+ \left(\frac{\eta_2}{\eta_0} \right)^{2t} r_2^2 + \left(\frac{\eta_1}{\eta_0}\right)^{2t}r_1^2\right)^2} dx_1\wedge d\bar{x}_1.
\]\\
\begin{prop}
\label{db}
$C_t$ is symplectic for $t\in [0,1]$. Namely, $C_0$ is symplectic isotropic to $C_1$ via the family $\{C_t\}_{t\in [0,1]}$ of smooth symplectic curves. More precisely, $(\tilde{\cal F}_t^* \omega_{\rm FS})|_{C_0}$ is smooth and is an $O(\epsilon)$-perturbation of $({\cal F}_t^* \omega_{\rm FS})|_{C_0}$.
\end{prop}
{\bf Proof:} According to lemma \ref{da} and proposition \ref{dm}, it is sufficient to show that $(\tilde{\cal F}_t^* \omega_{\rm FS})|_{C_0}$ is an $O(\epsilon)$-perturbation of $({\cal F}_t^* \omega_{\rm FS})|_{C_0}$.\\

Since $(\tilde{\cal F}_t^* \omega_{\rm FS})|_{C_0}$ and $({\cal F}_t^* \omega_{\rm FS})|_{C_0}$ coincide away from an $\epsilon$-neighborhood of ${\rm Sing}(C_0)$, with the help of symmetry, the cases that remain to be verified are $\epsilon$-neighborhoods of $\{r_1=r_2\leq 1-\epsilon\}$, $\{r_2=1,0\leq r_1\leq 1-\epsilon\}$ and $\{r_1=r_2=1\}$. On this neighborhoods, it is easy to observe that $\eta_1 = 1 + O(\epsilon)$, $\eta_2 = 1 + O(\epsilon)$, $\eta_1 = r_2 + O(\epsilon)$. Compare the expressions of $\tilde{\omega}_t$ and $({\cal F}_t^* \omega_{\rm FS})|_{C_0}$, we have that $\tilde{\omega}_t$ is an $O(\epsilon)$-perturbation of $({\cal F}_t^* \omega_{\rm FS})|_{C_0}$. Only thing remains to be shown is $R_t = O(\epsilon)$.\\

In an $\epsilon$-neighborhood of $\{r_1=r_2\leq 1-\epsilon\}$, $\eta_k = 1 + O(\epsilon)$ for $k=1,2$, $\lambda_1=\lambda_2=0$ and ${\rm Re}\left(\frac{x_2}{x_1}\right) - {\rm Re}\left(\frac{x_1}{x_2}\right) = O(\epsilon)$. Consequently, $R_t = t(1 -\lambda_0) R_{t,0} = O(\epsilon)$.\\

In an $\epsilon$-neighborhood of $\{r_2=1,0\leq r_1\leq 1-\epsilon\}$, $\eta_k = 1 + O(\epsilon)$ for $0\leq k \leq 2$, $\lambda_1=0$, $1-\lambda_0=0$ and ${\rm Re}(x_1) - {\rm Re} \left(\frac{x_1}{x_2}\right) = O(\epsilon)$. Consequently, $R_t = t\lambda_2 R_{t,2} = O(\epsilon)$.\\

In an $\epsilon$-neighborhood of $\{r_1=r_2=1\}$, $\eta_k = 1 + O(\epsilon)$ for $0\leq k \leq 2$, ${\rm Re}\left(\frac{x_2}{x_1}\right) - {\rm Re}\left(\frac{x_1}{x_2}\right) = O(\epsilon)$, ${\rm Re}(x_1) - {\rm Re} \left(\frac{x_1}{x_2}\right) = O(\epsilon)$, ${\rm Re}(x_2) - {\rm Re} \left(\frac{x_2}{x_1}\right) = O(\epsilon)$. Consequently, $R_{t,k} = O(\epsilon)$ for $0\leq k \leq 2$ and $R_t = O(\epsilon)$.
\hfill$\Box$\\
\begin{theorem}
\label{dn}
There exists a family of Hamiltonian diffeomorphism $H_t: {\mathbb{CP}^2}\rightarrow{\mathbb{CP}^2}$ such that $H_t(C_0)=C_t$ and $H_t$ is identity away from an arbitrary small neighborhood of $C_{[0,1]}$. The perturbed moment map (Lagrangian fibration) $\hat{F} = F\circ H_1$ is smooth and satisfies $\hat{F}(C_0) = \Gamma$ (the ``Y" shaped graph with a 3-valent vertex $v_0$) away from a small neighborhood of $v_0$. ($H_t$ can be made to be identity on $\partial \mathbb{CP}^2$ with the expense of smoothness of $\hat{F}$ at $\partial C_0 := \partial \mathbb{CP}^2 \cap C_0$.)
\end{theorem}
{\bf Proof:} Proposition \ref{db} implies that $C_0$ is smoothly symplectic isotropic to $C_1$ via the family $\{C_t\}_{t\in [0,1]}$. By the extension theorem (theorem 6.1) in \cite{lag2}, we can get a family of $C^\infty$ Hamiltonian diffeomorphism $H_t: {\mathbb{CP}^2}\rightarrow{\mathbb{CP}^2}$ such that $H_t(C_0)=C_t$ and $H_t$ is identity away from an arbitrary small neighborhood of $C_{[0,1]}$. To ensure that $H_t$ leaves $\partial \mathbb{CP}^2$ invariant, we need to use the extension theorem (theorem 6.6) in \cite{lag2}. Then $H_t$ can only be made $C^\infty$ away from the three intersection points of $C_t$ and $\partial \mathbb{CP}^2$.
\hfill$\Box$\\

\subsection{The optimal smoothness}
$\hat{F}$ constructed in section 4.2 is smooth. ($\hat{F}$ is not smooth at $\partial C_0  = \partial \mathbb{CP}^2 \cap C_0$ if $\hat{F}$ is required to be equal to $F$ on $\partial \mathbb{CP}^2$. This non-smoothness is due to the fact that $C_0$ is not symplectically normal crossing to $\partial \mathbb{CP}^2$ under $\omega_{\rm FS}$ and can be cured by modifying $\omega_{\rm FS}$ near $\partial C_0$ so that $C_0$ is symplectically normal crossing to $\partial \mathbb{CP}^2$.) The trade off is that $\hat{F}(C_0) = \Gamma$ (the ``Y" shaped graph with a 3-valent vertex $v_0$) away from a small neighborhood of $v_0$.\\

$\hat{F}$ constructed in section 4.1 satisfies $\hat{F}(C_0) = \Gamma$, but is only piecewise smooth and is not smooth at ${\rm Sing}(C_0)$. A natural question is: What is the optimal smoothness that $\hat{F}$ can achieve if we insist $\hat{F}(C_0) = \Gamma$? Clearly,  $\hat{F}$ can not be smooth over $v_0$. In this section, we will show that $\hat{F}$ can be made smooth over $\Gamma$ away from $v_0$. More precisely, let ${\rm Sing}_o(C_0) = \hat{F}^{-1}(v_0) \cap {\rm Sing}(C_0)$, we will show that $\hat{F}$ can be made smooth away from ${\rm Sing}_o(C_0)$. ($\hat{F}$ is not smooth at $\partial C_0$, if $\partial \mathbb{CP}^2$ is required to be fixed under $\hat{F}$.)\\

Let $b(a)$ be a smooth non-decreasing function satisfying $b(a) =0$ for $a\leq 0$, $b(a) >0$ for $a> 0$, $b(a)=1$ for $a\geq \sqrt{\epsilon}$ and $b'(a) \leq C/\sqrt{\epsilon}$. We may modify the definition of $C_t$ to consider $C_t = \tilde{\cal F}_t(C_0)$, where

\[
\tilde{\cal F}_t (x_1,x_2) = \left(\textstyle \left(\frac{\eta_1}{\eta_0}\right)^tx_1, \left(\frac{\eta_2}{\eta_0}\right)^tx_2\right),
\]
\[
\textstyle \log \eta_2 = \log r_1 h\left(\frac{\log r_1}{b_1}\right), \ b_1 = b(\log (r_1/r_2^2)),
\]
\[
\textstyle \log \eta_1 = \log r_2 h\left(\frac{\log r_2}{b_2}\right), \ b_2 = b(\log (r_2/r_1^2)),
\]
\[
\textstyle \log \eta_0 = \log r_1 h\left(\frac{\log (r_1/r_2)}{b_0}\right) + \log r_2 h\left(\frac{\log (r_2/r_1)}{b_0}\right), \ b_0 = b(\log (r_1r_2)).
\]\\
Notice that $\tilde{\cal F}_t$ here coincides with $\tilde{\cal F}_t$ in section 4.2 away from a $\sqrt{\epsilon}$-neighborhood of ${\rm Sing}_o(C_0)$, coincides with ${\cal F}_t$ in section 4.1 near ${\rm Sing}_o(C_0)$ away from a $\sqrt{\epsilon}$-neighborhood of $v_0$. Therefore, the only new construction of $\tilde{\cal F}_t$ is over a $\sqrt{\epsilon}$-neighborhood of $v_0$.\\

Assume $\lambda_0 = \lambda\left(\frac{\log r_2 - \log r_1}{b_0}\right)$, $\lambda_1 = \lambda\left(\frac{\log r_1}{b_1}\right)$, $\lambda_2 = \lambda\left(\frac{\log r_2}{b_2}\right)$. Then\\
\[
\textstyle \frac{d\eta_2}{\eta_2} = \lambda_1\frac{dr_1}{r_1} - \beta_1, \ \frac{d\eta_1}{\eta_1}= \lambda_2\frac{dr_2}{r_2} - \beta_2,
\]
\[
\textstyle \frac{d\eta_0}{\eta_0} = \frac{dr_1}{r_1} + \lambda_0\left(\frac{dr_2}{r_2} - \frac{dr_1}{r_1}\right) - \beta_0.
\]
\begin{lm}
\label{du}
$\beta_i = O(\sqrt{\epsilon})$ for $i = 1,2,3$.
\end{lm}
{\bf Proof:}
\[
\textstyle \beta_1 = \left[\frac{\log r_1}{b_1} h'\left(\frac{\log r_1}{b_1}\right)\right] \left[ \frac{\log r_1}{b_1}b'\left(\log \frac{r_1}{r_2^2}\right)\right] \left(\frac{dr_1}{r_1} - 2\frac{dr_2}{r_2}\right).
\]

Notice that $h'\left(\frac{\log r_1}{b_1}\right) \not=0$ only when $\frac{\log r_1}{b_1} \leq \epsilon$. Hence

\[
\textstyle \left[\frac{\log r_1}{b_1} h'\left(\frac{\log r_1}{b_1}\right)\right] = O(1),\  \left[ \frac{\log r_1}{b_1}b'\left(\log \frac{r_1}{r_2^2}\right)\right] = O(\sqrt{\epsilon}).
\]

Consequently, $\beta_1 = O(\sqrt{\epsilon})$. The verifications for $\beta_2$ and $\beta_3$ are similar.
\hfill$\Box$\\

By similar computation as in section 4.2, we get\\
\[
\frac{(\tilde{\cal F}_t^* \omega_{\rm FS})|_{C_0}}{dx_1\wedge d\bar{x}_1} = \frac{\tilde{\omega}_t}{dx_1\wedge d\bar{x}_1} + tR_t + tB_t \geq \frac{1}{6} + O(\sqrt{\epsilon}),
\]

where $B_t$ is linear on $\{\beta_i\}_{i=1}^3$ and $B_t = O(\sqrt{\epsilon})$.\\
\begin{prop}
\label{ds}
$C_t$ is symplectic for $t\in [0,1]$. Namely, $C_0$ is symplectic isotropic to $C_1$ via the family $\{C_t\}_{t\in [0,1]}$ of smooth symplectic curves. More precisely, $(\tilde{\cal F}_t^* \omega_{\rm FS})|_{C_0}$ is smooth away from ${\rm Sing}_o(C_0)$ and is an $O(\sqrt{\epsilon})$-perturbation of $({\cal F}_t^* \omega_{\rm FS})|_{C_0}$.
\end{prop}
{\bf Proof:} This proposition is a direct consequence of the above computation, lemma \ref{da}, propositions \ref{dm} and \ref{db} together with the additional estimate $B_t = O(\sqrt{\epsilon})$ implied by lemma \ref{du}.
\hfill$\Box$\\
\begin{theorem}
\label{dt}
There exists a family of Hamiltonian diffeomorphism $H_t: {\mathbb{CP}^2} \rightarrow {\mathbb{CP}^2}$ such that $H_t(C_0)=C_t$ and $H_t$ is identity away from an arbitrary small neighborhood of $C_{[0,1]}$. $\hat{F} = F\circ H_1$ satisfies $\hat{F}(C_0) = \Gamma$ (the ``Y" shaped graph with a 3-valent vertex) and is smooth away from ${\rm Sing}_o(C_0)$. ($H_t$ can be made to be identity on $\partial \mathbb{CP}^2$ with the expense of smoothness of $\hat{F}$ at $\partial \mathbb{CP}^2 \cap C_0$.)
\end{theorem}
{\bf Proof:} The proof is essentially the same as the proofs of theorems \ref{dc} except here $C_t$ is decomposed into 3 (instead of 6) smooth symmetric pieces, lemma \ref{da} and proposition \ref{dm} is replaced by proposition \ref{ds} and ${\rm Sing}(C_0)$ is replaced by ${\rm Sing}_o(C_0)$.
\hfill$\Box$\\

\se{String diagram and Feynman diagram}
In this section, we will naturally combine the localization technique of section 3, which reduces the curves (string diagram) locally to individual pair of pants, with the explicit perturbation technique of section 4 to perturb the moment map $F_{\delta^w}$, so that the perturbed moment map will map $C_{s^{\delta^w}}$ to a graph. This is a very interesting analogue of the relation of string diagrams in string theory and Feynman diagrams in quantum mechanics.\\

In general, given a simplicial decomposition $Z\in \tilde{Z}$ of $\Delta$, take a weight $w\in \tau^0_Z$, according to proposition \ref{cf}, we have

\[
P_\Sigma = \bigcup_{S\in Z} U^S_\epsilon.
\]

According to results in \cite{lag2}, the perturbation of the moment map can be reduced to the perturbation of the pair $(C_{s^{\delta^w}}, \omega_{\delta^w})$ of symplectic curve and symplectic form. For each $S \in Z^{\rm top}$, locally in $U^S_\epsilon$, $(C_{s^{\delta^w}} \cap U^S_\epsilon, \omega_{\delta^w}|_{U^S_\epsilon})$ is a close approximation of the line and the Fubini-Study \k form discussed in section 4. Namely, the construction in section 4 can be viewed as local model for construction here. In the following, we will start with some modification of the local model in section 4, then we will apply the modified local model to perturb $C_{s^{\delta^w}}$. For such purpose, $\omega_{\delta^w}$ also need to be perturbed suitably.\\

\subsection{Modified local models}
Consider a smooth non-negative non-decreasing function $\gamma_\epsilon(u)$, such that $\gamma_\epsilon(u) = 0$ for $\sqrt{u}\leq A_1\epsilon$ and $\gamma_\epsilon(u)=1$ for $\sqrt{u}\geq A_2\epsilon$. $A_1,A_2$ are positive constants satisfying $1 < A_1 <A_2 <|\Delta|$. Let $\gamma_{\epsilon,t} (u) = t\gamma_\epsilon (u) + (1-t)$ and\\
\[
\eta_1 = \max(1,r_2),\ \ \eta_2 = \max(1,r_1),\ \ \eta_0 = \max(r_1,r_2).
\]
\begin{prop}
\label{dj}
$C_t = p_t^{-1}(0)$ is symplectic curve under the Fubini-Study \k form for $t\in [0,1]$, where

\[
\textstyle p_t(x) = \gamma_{\epsilon,t} \left(\frac{r_1^2}{\eta_1^2}\right)x_1 + \gamma_{\epsilon,t} \left(\frac{r_2^2}{\eta_2^2}\right)x_2 + \gamma_{\epsilon,t} \left(\frac{1}{\eta_0^2}\right)=0.
\]

Namely, the family $\{C_t\}_{t\in [0,1]}$ is a symplectic isotopy from $C_0 = \{(x_1,x_2): x_1+x_2+1=0\}$ to

\[
C_1 = \left\{\textstyle \left(x_1, x_2\right): \gamma_\epsilon\left(\frac{r_1^2}{\eta_1^2}\right)x_1 + \gamma_\epsilon\left(\frac{r_2^2}{\eta_2^2}\right)x_2 +\gamma_\epsilon\left(\frac{1}{\eta_0^2}\right)=0\right\}.
\]
\end{prop}
{\bf Proof:} By symmetry, we only need to verify that $C_t$ is symplectic in the region $1\geq |x_2|\geq |x_1|$, where

\[
p_t(x) = \gamma_{\epsilon,t} (|x_1|^2)x_1 + x_2 + 1 =0.
\]

Since $C_t$ is a complex curve away from the region $\{A_1\epsilon \leq |x_1| \leq A_2\epsilon\}$, we only need to verify that $C_t \cap \{A_1\epsilon \leq |x_1| \leq A_2\epsilon\}$ is symplectic.\\

Recall the \k form of the Fubini-Study metric is

\[
\omega_{\rm FS} = \frac{ dx_1\wedge d\bar{x}_1 + dx_2\wedge d\bar{x}_2 + (x_2 dx_1 - x_1 dx_2)\wedge (\bar{x}_2d\bar{x}_1-\bar{x}_1d\bar{x}_2)}{(1+|x|^2)^2}.
\]

When restricted to $C_t \cap \{A_1\epsilon \leq |x_1| \leq A_2\epsilon\}$,

\[
\omega_{\rm FS} = \frac{1}{2}dx_1\wedge d\bar{x}_1 + \frac{1}{4} dx_2\wedge d\bar{x}_2 + O(\epsilon)
\]
\[
= \left(\frac{1}{2} + \frac{1}{4} [\gamma_{\epsilon,t}(|x_1|^2)^2 + \gamma_{\epsilon,t}(|x_1|^2)t\gamma_\epsilon^{\#} (|x_1|^2)]\right) dx_1\wedge d\bar{x}_1 + O(\epsilon)
\]
\[
\geq \frac{2 + (1-t)^2}{4} dx_1\wedge d\bar{x}_1 + O(\epsilon),
\]

where $\gamma_\epsilon^{\#}(|x_1|^2) = 2|x_1|^2\gamma_\epsilon'(|x_1|^2)$. Therefore $C_t$ is symplectic.
\hfill$\Box$\\
\begin{prop}
\label{dh}
$C_t = {\cal F}_t (C_0)$ is symplectic for $t\in [0,1]$, where
\[
\textstyle {\cal F}_t (x_1,x_2) = \left(\left(\frac{\eta_1}{\eta_0}\right)^tx_1, \left(\frac{\eta_2}{\eta_0}\right)^tx_2\right),
\]
\begin{equation}
\label{ea}
\textstyle C_0 = \left\{(x_1, x_2): \gamma_\epsilon\left(\frac{r_1^2}{\eta_1^2}\right)x_1 + \gamma_\epsilon\left(\frac{r_2^2}{\eta_2^2}\right)x_2 +\gamma_\epsilon \left(\frac{1}{\eta_0^2}\right) =0 \right\}.
\end{equation}
\end{prop}
{\bf Proof:} By symmetry, we only need to verify that $C_t$ is symplectic in the region $1\geq |x_2|\geq |x_1|$, which is one of the six symmetric regions that together form ${\mathbb{CP}^2}$. In the region $1\geq |x_2|\geq |x_1|$,

\[
C_t = \left\{\textstyle \left( \left(\frac{1}{r_2}\right)^tx_1,\left(\frac{1}{r_2}\right)^tx_2\right): \gamma_\epsilon(|x_1|^2)x_1 + x_2 +1=0\right\}
\]

$\gamma_\epsilon(|x_1|^2)x_1 + x_2 +1=0$ implies that

\[
dx_2 = - \gamma_\epsilon dx_1 - x_1 d\gamma_\epsilon.
\]

Hence\\
\[
\textstyle \frac{dr_2}{r_2}= {\rm Re}\left(\frac{dx_2}{x_2}\right)= -\gamma_\epsilon{\rm Re}\left(\left(\frac{x_1}{x_2}\right)\frac{dx_1}{x_1}\right) - {\rm Re}\left(\frac{x_1}{x_2}\right)d\gamma_\epsilon.
\]\\
\[
\textstyle d\left(\left(\frac{1}{r_2}\right)^tx_1\right)\wedge d\left(\left(\frac{1}{r_2}\right)^t\bar{x}_1\right) = \left(\frac{1}{r_2}\right)^{2t}\left(dx_1\wedge d\bar{x}_1 +t(x_1d\bar{x}_1 - \bar{x}_1dx_1)\wedge \frac{dr_2}{r_2}\right)
\]
\[
\textstyle =\left(\frac{1}{r_2}\right)^{2t} \left(1 + t(\gamma_\epsilon + \gamma_\epsilon^{\#}){\rm Re}\left(\frac{x_1}{x_2}\right)\right)dx_1\wedge d\bar{x}_1,
\]\\
\[
\textstyle d\left(\left(\frac{1}{r_2}\right)^tx_2\right)\wedge d\left(\left(\frac{1}{r_2}\right)^t\bar{x}_2\right) = \left(\frac{1}{r_2}\right)^{2t}\left(dx_2\wedge d\bar{x}_2 +t(x_2d\bar{x}_2 - \bar{x}_2dx_2)\wedge \frac{dr_2}{r_2}\right)
\]
\[
\textstyle =\left(\frac{1}{r_2}\right)^{2t} (1 - t)dx_2\wedge d\bar{x}_2 = \left(\frac{1}{r_2}\right)^{2t} (1 - t)(\gamma_\epsilon^2 + \gamma_\epsilon\gamma_\epsilon^{\#})dx_1\wedge d\bar{x}_1 \geq 0,
\]\\
\[
\textstyle \left( \left(\frac{1}{r_2}\right)^tx_2\right)d\left( \left(\frac{1}{r_2}\right)^tx_1\right) - \left( \left(\frac{1}{r_2}\right)^tx_1\right)d\left( \left(\frac{1}{r_2}\right)^tx_2\right)
\]
\[
\textstyle = \left(\frac{1}{r_2}\right)^{2t}(x_2dx_1 - x_1dx_2) = -\left(\frac{1}{r_2}\right)^{2t} (dx_1 + x_1^2d\gamma_\epsilon).
\]\\
By restriction to $C_t$ we get

\[
\frac{\omega_{\rm FS}|_{C_t}}{dx_1\wedge d\bar{x}_1} \geq \frac{\left(\frac{1}{r_2}\right)^{2t}\left(1 + t(\gamma_\epsilon + \gamma_\epsilon^{\#}){\rm Re}\left(\frac{x_1}{x_2}\right)\right) + \left(\frac{1}{r_2}\right)^{4t}(1 - {\rm Re}(x_1)\gamma_\epsilon^{\#})}{\left(1+r_2^{2-2t} + \left(\frac{r_1}{r_2}\right)^{2t} r_1^{2-2t}\right)^2}
\geq \frac{1}{6} + O(\epsilon).
\]

The reason is that $0\leq \gamma_\epsilon \leq 1$, $\gamma_\epsilon^{\#}{\rm Re}\left(\frac{x_1}{x_2}\right) = O(\epsilon)$ and ${\rm Re}(x_1)\gamma_\epsilon^{\#} = O(\epsilon)$. Therefore $C_t$ is symplectic.
\hfill$\Box$\\

{\bf Remark:} Let $U^{\mathbb{CP}^2}_\epsilon = \{r_1 \leq \epsilon\eta_1, r_2 \leq \epsilon\eta_2, 1 \leq \epsilon\eta_0\} \subset \mathbb{CP}^2$. It is easy to observe that outside of $U^{\mathbb{CP}^2}_\epsilon$, $C_t$ in proposition \ref{dh} is equal to $\{x_2 + 1 =0\}$ when $|x_1|$ is small, equal to $\{x_1 + 1 =0\}$ when $|x_2|$ is small, equal to $\{x_1 + x_2 =0\}$ when $|x_1|,|x_2|$ are large. Namely, $C_t$ outside of $U^{\mathbb{CP}^2}_\epsilon$ is toric, $F(C_t \cap (\mathbb{CP}^2 \setminus U^{\mathbb{CP}^2}_\epsilon))$ is 1-dimensional, independent of $t$ and is the union of the 3 end segments of the ``Y" shaped graph. Also the image of $C_1$ under any moment map is a 1-dimensional graph of ``Y" shape.\\

The following is the analogue of theorem \ref{dc} for our modified local model.\\
\begin{theorem}
\label{dq}
There exists a family of piecewise smooth Lipschitz Hamiltonian diffeomorphism $H_t: {\mathbb{CP}^2}\rightarrow{\mathbb{CP}^2}$ such that $H_t(C_0)=C_t$ and $H_t$ is identity away from an arbitrary small neighborhood of $C_{[0,1]}$ or away from $U^{\mathbb{CP}^2}_\epsilon$. The perturbed moment map (Lagrangian fibration) $\hat{F} = F\circ H_1$ satisfies $\hat{F}(C_0) = \Gamma$ (the ``Y" shaped graph with a 3-valent vertex).
\end{theorem}
{\bf Proof:} The proof is essentially the same as the proof of theorem \ref{dc} except for the proof of $H_t$ being the identity map when restricted to $\mathbb{CP}^2 \setminus U^{\mathbb{CP}^2}_\epsilon$, which is based on the fact that ${\cal F}_t$ restricts to identity map on $C_0 \setminus U^{\mathbb{CP}^2}_\epsilon$.
\hfill$\Box$\\

To deal with the cases of smooth and optimal smoothness discussed in sections 4.2 and 4.3, we may take $C_t = \tilde{\cal F}_t(C_0)$, where we take $C_0$ in (\ref{ea}) and $\tilde{\cal F}_t$ in either section 4.2 or section 4.3. (Notice that in the region where $C_0$ is modified, $\tilde{\cal F}_t$ in sections 4.2 and 4.3 coincide.)\\
\begin{prop}
\label{eb}
$C_t = \tilde{\cal F}_t (C_0)$ is symplectic for $t\in [0,1]$, where $C_0$ is defined in (\ref{ea}).
\end{prop}
{\bf Proof:} By symmetry, we only need to verify that $C_t$ is symplectic in the region, where $|x_1| \leq |x_2| \leq 1$ and $\gamma_\epsilon (|x_1|^2)<1$. In this region, we have $|x_1| = O(\epsilon)$ and  $x_2 = -1 + O(\epsilon)$. Consequently, $\lambda_0 -1 = \lambda_1 = 0$, $\eta_2 =1$, $\eta_1 = 1 + O(\epsilon)$ and $\eta_0 = r_2 = 1 + O(\epsilon)$.\\
\[
\textstyle \frac{d\eta_2}{\eta_2} = 0, \frac{d\eta_1}{\eta_1}= \lambda_2\frac{dr_2}{r_2}, \frac{d\eta_0}{\eta_0}= \frac{dr_2}{r_2}.
\]

$\gamma_\epsilon(|x_1|^2)x_1 + x_2 +1=0$ implies that

\[
dx_2 = - \gamma_\epsilon dx_1 - x_1 d\gamma_\epsilon.
\]

Hence\\
\[
\textstyle \frac{dr_2}{r_2}= {\rm Re}\left(\frac{dx_2}{x_2}\right)= -\gamma_\epsilon{\rm Re}\left(\left(\frac{x_1}{x_2}\right)\frac{dx_1}{x_1}\right) - {\rm Re}\left(\frac{x_1}{x_2}\right)d\gamma_\epsilon.
\]
\[
dx_2\wedge d\bar{x}_2 = \gamma_\epsilon(\gamma_\epsilon + \gamma_\epsilon^{\#})dx_1\wedge d\bar{x}_1
\]\\
\[
\textstyle d\left(\left(\frac{\eta_1}{\eta_0}\right)^tx_1\right) = \left(\frac{\eta_1}{\eta_0}\right)^t\left(dx_1 + tx_1\left(\frac{d\eta_1}{\eta_1} - \frac{d\eta_0}{\eta_0}\right)\right).
\]
\[
\textstyle \frac{d\eta_1}{\eta_1} - \frac{d\eta_0}{\eta_0} = - (1 - \lambda_2)\frac{dr_2}{r_2}.
\]
\begin{eqnarray*}
&&\textstyle d\left(\left(\frac{\eta_1}{\eta_0}\right)^tx_1\right)\wedge d\left(\left(\frac{\eta_1}{\eta_0}\right)^t\bar{x}_1\right)\\
&=& \textstyle \left(\frac{\eta_1}{\eta_0}\right)^{2t}\left(dx_1\wedge d\bar{x}_1 -t(x_1d\bar{x}_1 - \bar{x}_1dx_1)\wedge \left(\frac{d\eta_1}{\eta_1} - \frac{d\eta_0}{\eta_0}\right)\right)\\
&=&\textstyle \left(\frac{\eta_1}{\eta_0}\right)^{2t} \left(1 + (\gamma_\epsilon + \gamma_\epsilon^{\#})(1 - \lambda_2)t{\rm Re}\left(\frac{x_1}{x_2}\right)\right)dx_1\wedge d\bar{x}_1.
\end{eqnarray*}

\[
\textstyle d\left(\left(\frac{\eta_2}{\eta_0}\right)^tx_2\right) = \left(\frac{\eta_2}{\eta_0}\right)^t\left(dx_2 + tx_2\left(\frac{d\eta_2}{\eta_2} - \frac{d\eta_0}{\eta_0}\right)\right).
\]
\[
\textstyle \frac{d\eta_2}{\eta_2} - \frac{d\eta_0}{\eta_0} = - \frac{dr_2}{r_2}.
\]
\begin{eqnarray*}
&&\textstyle d\left(\left(\frac{\eta_2}{\eta_0}\right)^tx_2\right)\wedge d\left(\left(\frac{\eta_2}{\eta_0}\right)^t\bar{x}_2\right)\\
&=& \textstyle \left(\frac{\eta_2}{\eta_0}\right)^{2t}\left(dx_2\wedge d\bar{x}_2 +t(x_2d\bar{x}_2 - \bar{x}_2dx_2)\wedge \left(\frac{d\eta_2}{\eta_2} - \frac{d\eta_0}{\eta_0}\right)\right)\\
&=&\textstyle \left(\frac{\eta_2}{\eta_0}\right)^{2t} \gamma_\epsilon(\gamma_\epsilon + \gamma_\epsilon^{\#})(1 - t) dx_1\wedge d\bar{x}_1 \geq 0.
\end{eqnarray*}

\[
\textstyle \alpha = \left( \left(\frac{\eta_2}{\eta_0}\right)^tx_2\right)d\left( \left(\frac{\eta_1}{\eta_0}\right)^tx_1\right) - \left( \left(\frac{\eta_1}{\eta_0}\right)^tx_1\right)d\left( \left(\frac{\eta_2}{\eta_0}\right)^tx_2\right)
\]
\[
\textstyle = \left(\frac{\eta_2\eta_1}{\eta_0^2}\right)^{t}\left(x_2dx_1 - x_1dx_2 + tx_1x_2\left(\frac{d\eta_1}{\eta_1} - \frac{d\eta_2}{\eta_2}\right)\right)
\]
\[
\textstyle = \left(\frac{\eta_2\eta_1}{\eta_0^2}\right)^{t} \left(-dx_1 + x_1^2 d\gamma_\epsilon + tx_1x_2 \lambda_2\frac{dr_2}{r_2}\right) = -\left(\frac{\eta_2\eta_1}{\eta_0^2}\right)^{t} dx_1 + O(|x_1|).
\]\\
\[
\textstyle \alpha \wedge \bar{\alpha} = \left(\frac{\eta_2\eta_1}{\eta_0^2}\right)^{2t} dx_1d\bar{x}_1 + O(|x_1|).
\]

By restriction to $C_t$ we get\\

\begin{tabular}{cl}
\centerline{$\displaystyle \frac{(\tilde{\cal F}_t^* \omega_{\rm FS})|_{C_0}}{dx_1\wedge d\bar{x}_1} \geq \frac{\left(\frac{\eta_1}{\eta_0}\right)^{2t} + \left(\frac{\eta_2\eta_1}{\eta_0^2}\right)^{2t} + O(|x_1|)}{\left(1+\left(\frac{\eta_2}{\eta_0}\right)^{2t}r_2^2 + \left(\frac{\eta_1}{\eta_0}\right)^{2t}r_1^2\right)^2} \geq \frac{1}{2} + O(\epsilon).$} &  \hspace*{-.4in} $\Box$
\end{tabular}\\\\

For $C_t = \tilde{\cal F} (C_0)$, where $\tilde{\cal F}$ is taken from section 4.3, we have\\
\begin{theorem}
\label{ec}
$\hat{F}$ in theorem \ref{dq} can be made smooth away from ${\rm Sing}_o(C_0)$.
\end{theorem}
{\bf Proof:} The proof is essentially the same as the proof of theorem \ref{dc} except that $C_0$ is decomposed into 3 pieces with boundaries in ${\rm Sing}_o(C_0)$. The proof of $H_t$ being the identity map when restricted to $\mathbb{CP}^2 \setminus U^{\mathbb{CP}^2}_\epsilon$ is based on the fact that $\tilde{\cal F}_t$ restricts to identity map on $C_0 \setminus U^{\mathbb{CP}^2}_\epsilon$.
\hfill$\Box$\\

{\bf Remark:} There is also a version of theorem \ref{ec} as analogue of theorem \ref{dn} when $\tilde{\cal F}$ is taken from section 4.2.\\

\subsection{Perturbation of symplectic curve and form}
For $m\in \Delta$, let

\[
\Delta_m = \{m' \in \Delta| \{m,m'\}\in Z\}.
\]

Choose $\check{\epsilon}$ such that $\delta^a \leq \check{\epsilon} \leq \epsilon$. Define

\[
\hat{s}_m = \gamma_\epsilon (\rho_m)s_m,\ \check{s}_m = \left[ 1 - \gamma_{\check{\epsilon}} \left(\max_{m'\not\in \Delta_m}\left(\rho_{m'} \right)\right)\right] s_m,
\]
\[
\hat{s}^{\delta^w} = \sum_{m\in \Delta} \delta^{w_m}a_m\hat{s}_m,\ \check{s}^{\delta^w} = \sum_{m\in \Delta} \delta^{w_m}a_m\check{s}_m.
\]
\[
\check{\omega}_{\delta^w} = \partial \bar{\partial} \check{h}_{\delta^w}, \mbox{ where } \check{h}_{\delta^w} = \log|\check{s}^{\delta^w}|_\Delta^2, \ |\check{s}^{\delta^w}|_\Delta^2 = \sum_{m\in \Delta} |\delta^{w_m}\check{s}_m|_\Delta^2.
\]

\begin{prop}
\label{do}
$\check{\omega}_{\delta^w}$ is a \k form on $P_\Sigma$ near $C_t = s_t^{-1}(0)$ for $t\in [0,1]$, where $s_t = t\hat{s}_{\delta^w} + (1-t)s_{\delta^w}$.
\end{prop}
{\bf Proof:} For $x\in P_\Sigma$, let $\rho_{m_i}(x)$ for a $m_i\in \Delta$ be the $i$-th largest among $\{\rho_m(x)\}_{m\in \Delta}$. Since $S_x$ is non-empty, we have $m_1 \in S_x$ and $\rho_{m_1}(x) \geq 1/|\Delta| - \epsilon$. If $x \in C_t$, it is easy to derive from the equation of $C_t$ that $\rho_{m_2}(x) \geq 1/|\Delta|^2 - \epsilon/|\Delta|$ and $m_2\in S_x$ when $\epsilon$ is small.\\

If $\{m_1,m_2\} \not\subset \Delta_m$, then $\displaystyle \max_{m'\not\in \Delta_m}\left(\rho_{m'} (x) \right) \geq \rho_{m_2} (x) > |\Delta| \check{\epsilon}$ when $\check{\epsilon}$ is small. Hence $\check{s}_m (x) = 0$.\\

If $\{m_1,m_2\} \subset \Delta_m$ and $\check{s}_m \not= s_m$, then there exists $m' \not\in \Delta_m$ such that $\rho_{m'} > \check{\epsilon}$. Hence $\check{S}_x = \{m_1,m_2,m'\}$, $\check{s}_{m'} = s_{m'}$ and $m_3 = m'$, where $\check{S}_x$ is defined as $S_x$ with $\epsilon$ replaced by $\check{\epsilon}$. Consequently, $\rho_m(x) = O(\delta^a)$ and $\check{\omega}_{\delta^w}(x)$ is an $O(\delta^a /\check{\epsilon})$-perturbation of $\check{\omega}_{\delta^w}^{\check{S}_x}(x)$. When $\delta^a /\check{\epsilon}$ is small, $\check{\omega}_{\delta^w}$ is a \k form at $x$.\\

The remaining case is when $\check{s}_m = s_m$ for $m\in S' = \{m_1,m_2, m', m''\}$ and $\check{s}_m = 0$ for $m\not\in S'$, where $\{m', m''\}$ is uniquely determined by the relation $\{m_1,m_2\} \subset \Delta_{m'} \cap \Delta_{m''}$. Then $\check{\omega}_{\delta^w} (x) = \omega^{S'}_{\delta^w} (x)$ is clearly \k. Therefore $\check{\omega}_{\delta^w}$ is a \k form on $P_\Sigma$ near $C_t$.
\hfill$\Box$\\
\begin{prop}
\label{dk}
$C_t$ is symplectic curve under the \k form $\omega_t$ for $t\in [0,1]$, where $\omega_t = t\check{\omega}_{\delta^w} + (1-t)\omega_{\delta^w}$. Namely, the family $\{C_t\}_{t\in [0,1]}$ is a symplectic isotopy from $C_0 = C_{s_{\delta^w}}$ to $C_1 = C_{\hat{s}_{\delta^w}}$. Further more, there exists smooth symplectomorphisms $H_1: (P_\Sigma, \omega_{\delta^w}) \rightarrow (P_\Sigma, \check{\omega}_{\delta^w})$ such that $H_1(C_{s_{\delta^w}}) = C_{\hat{s}_{\delta^w}}$. ($H_1$ can be made to be identity on $\partial P_\Sigma$ with the expense of smoothness of $H_1$ at $C_0 \cap \partial P_\Sigma$.)
\end{prop}
{\bf Proof:} Proposition \ref{do} implies that $\omega_t$ are \k forms on $P_\Sigma$ near $C_t$. It is easy to see that $s_t$ is holomorphic outside of the union of $U^S_\epsilon$ for $S\in Z^{\rm top}$, where $C_t$ is automatically symplectic.\\

For each $S = \{ m_0,m_1,m_2\} \in Z^{\rm top}$, $\{z_i = \delta^{w_{m_i}} a_{m_i} s_{m_i}\}_{i=0}^2$ defines an open embedding $U^S_\epsilon \hookrightarrow {\mathbb{CP}^2}$, where $[z_0,z_1,z_2]$ is the homogeneous coordinate of ${\mathbb{CP}^2}$. Using the inhomogeneous coordinates $(x_1,x_2)$ of ${\mathbb{CP}^2}$ on $U^S_\epsilon$, $\hat{s}_{\delta^w}$ reduces to $p_1$ in proposition \ref{dj} and $s_{\delta^w}$ reduces to $p_0 = x_1 + x_2 +1$ in proposition \ref{dj} up to $O(\delta^+)$ terms (lemma \ref{de}). Hence $C_t$ here coincides with $C_t$ in proposition \ref{dj} inside $U^S_\epsilon \subset {\mathbb{CP}^2}$. When $\delta$ is small, by proposition \ref{dj}, $C_t$ is symplectic in $U^S_\epsilon$ with respect to $\omega_{\rm FS}$. Since $\check{\omega}_{\delta^w} = \omega_{\rm FS}$ when restricted to $U^S_\epsilon$, $C_t$ is symplectic in $U^S_\epsilon$ with respect to $\check{\omega}_{\delta^w}$.\\

For the second part of the proposition, Apply theorems 6.1 and 6.2 from \cite{lag2} (which though are conveniently formulated for our application here, are essentially well known along the line of Moser's theorem) to the symplectic isotopic family $\{(C_t, \omega_t)\}_{t\in [0,1]}$, we can construct a smooth symplectomorphism $H_1: (P_\Sigma, \omega_{\delta^w}) \rightarrow (P_\Sigma, \check{\omega}_{\delta^w})$ such that $H_1(C_{s_{\delta^w}}) = C_{\hat{s}_{\delta^w}}$. To satisfy $H_1|_{\partial P_\Sigma} = {\rm Id}_{\partial P_\Sigma}$, it is necessary to apply theorems 6.3 and 6.4 from \cite{lag2} and $H_1$ is piecewise smooth, $C^{0,1}$ and is smooth away from $C_0 \cap \partial P_\Sigma$.
\hfill$\Box$\\

When $S \in Z$ is a 1-simplex, $\Gamma_S$ is just the baricenter of $S$. Let $s(\Gamma_Z)$ (resp. $e(\Gamma_Z)$) denote the union of $\Gamma_S$ for those 1-simplex $S \in Z$ that is not in $\partial \Delta$ (resp. is in $\partial \Delta$).\\
\begin{prop}
\label{dl}
For each $S\in Z^{\rm top}$, we may modify $C_{\hat{s}_{\delta^w}}$ in $U^S_\epsilon$ according to proposition \ref{dh}, and keep $C_{\hat{s}_{\delta^w}}$ unchanged outside of the union of such $U^S_\epsilon$. In such way, we can construct a family of symplectic curves $\{C_t\}_{t\in [0,1]}$ under the symplectic form $\check{\omega}_{\delta^w}$, such that $C_0 = C_{\hat{s}_{\delta^w}}$ and $F_{\delta^w}(C_1) = \Gamma$ is a graph that coincides with $\Gamma_Z$ away from an $\epsilon$-neighborhood of $s(\Gamma_Z)$ and is an $O(\epsilon)$-perturbation of $\Gamma_Z$.
\end{prop}
{\bf proof:} It is straightforward to verify that the deformation defined in the proposition match on overlaping regions. Through similar discussion as in the remark after proposition \ref{dh}, it is easy to observe that $C_t$ is toric outside of the union of $U^S_\epsilon$ for $S\in Z^{\rm top}$, hence the moment map image of $C_t$ in this region is 1-dimensional, independent of $t$ and is inside a small neighborhood of $s(\Gamma_Z) \cap e(\Gamma_Z)$. For each $S\in Z^{\rm top}$, in $U^S_\epsilon$, as in the proof of proposition \ref{dk}, we have coordinates $(x_1,x_2)$, which reduces $C_t$ here to $C_t \subset \mathbb{CP}^2$ in proposition \ref{dh}. Hence the image of $C_1 \cap U^S_\epsilon$ under the moment map coincides with part of $\Gamma_S \subset \Gamma_Z$ according to proposition \ref{dh}.
\hfill$\Box$\\
\begin{theorem}
\label{dd}
There exists a piecewise smooth Lagrangian fibration $\hat{F}$ as perturbation of the moment map $F_{\delta^w}$ such that $\hat{F}|_{\partial P_\Sigma} = F_{\delta^w}|_{\partial P_\Sigma}$ and $\hat{F}(C_{s^{\delta^w}}) = \Gamma$ is a graph that coincides with $\Gamma_Z$ away from an $\epsilon$-neighborhood of $s(\Gamma_Z)$ and is an $O(\epsilon)$-perturbation of $\Gamma_Z$.
\end{theorem}
{\bf Proof:} According to proposition \ref{dk}, we can construct a smooth symplectomorphism $H_1: (P_\Sigma, \omega_{\delta^w}) \rightarrow (P_\Sigma, \check{\omega}_{\delta^w})$ such that $H_1(C_{s_{\delta^w}}) = C_{\hat{s}_{\delta^w}}$. One can make $H_1|_{\partial P_\Sigma} = {\rm Id}_{\partial P_\Sigma}$ with the expense of smoothness of $H_t$ at $C_0 \cap \partial P_\Sigma$.\\

For the symplectic isotopic family $\{C_t\}_{t\in [0,1]}$ under the symplectic form $\check{\omega}_{\delta^w}$ in proposition \ref{dl}, we may define $H_2$ in $U^S_\epsilon$ for $S\in Z^{\rm top}$ to be the $H_1$ in theorem \ref{dq} and extend by identity map outside the union of $U^S_\epsilon$ for $S\in Z^{\rm top}$. Then $H_2: (P_\Sigma, \check{\omega}_{\delta^w}) \rightarrow (P_\Sigma, \check{\omega}_{\delta^w})$ is piecewise smooth and $C^{0,1}$ symplectomorphism satisfying $H_2|_{\partial P_\Sigma} = {\rm Id}_{\partial P_\Sigma}$, $H_2(C_{\hat{s}_{\delta^w}}) = C_1$ such that $\check{F}_{\delta^w} (C_1) = \Gamma$ is a graph that is an $\epsilon$-perturbation of the graph $\Gamma_Z$.\\

Let $H = H_2 \circ H_1$. Then $H|_{\partial P_\Sigma} = {\rm Id}_{\partial P_\Sigma}$ and $\hat{F} = F_{\delta^w}\circ H$ is the desired perturbation of $F_{\delta^w}$.
\hfill$\Box$\\

{\bf Remark:} Theorems \ref{dd} and \ref{ca} of this paper are needed for the proofs in \cite{lag3}.\\
\begin{prop}
\label{ed}
For each $S\in Z^{\rm top}$, we may modify $C_{\hat{s}_{\delta^w}}$ in $U^S_\epsilon$ according to proposition \ref{eb}, and keep $C_{\hat{s}_{\delta^w}}$ unchanged outside of the union of such $U^S_\epsilon$. In such way, we can construct a family of symplectic curves $\{C_t\}_{t\in [0,1]}$ under the symplectic form $\check{\omega}_{\delta^w}$, such that $C_0 = C_{\hat{s}_{\delta^w}}$ and $F_{\delta^w}(C_1) = \Gamma$ is a graph that coincides with $\Gamma_Z$ away from an $\epsilon$-neighborhood of $s(\Gamma_Z)$ and is an $O(\epsilon)$-perturbation of $\Gamma_Z$.
\end{prop}
{\bf proof:} The proof is the same as the proof of proposition \ref{dl} except that proposition \ref{dh} is replaced with proposition \ref{eb}.
\hfill$\Box$\\
\begin{theorem}
\label{ee}
$\hat{F}$ in theorem \ref{dd} can be made smooth away from $C_0 \cap \hat{F}^{-1} (v(\Gamma_Z))$ and $C_0 \cap \partial P_\Sigma$, where $v(\Gamma_Z)$ is the set of 3-valent vertices of $\Gamma_Z$.
\end{theorem}
{\bf Proof:} The proof is the same as the proof of theorem \ref{dd} except that proposition \ref{dl} (resp. theorem \ref{dq}) is replaced with proposition \ref{ed} (resp. theorem \ref{ec}).
\hfill$\Box$\\

{\bf Remark:} In this theorem, $\hat{F}$ achieved optimal smoothness possible. This result is a significant improvement over theorem \ref{dd}, and should play an important role in improving the Lagrangian torus fibration of quintic Calabi-Yau constructed in \cite{lag3} to optimal smoothness. We hope to come back to such improvement of \cite{lag3} in a future paper.\\

Theorems \ref{dd} and \ref{ee} concern the partial secondary fan, where $Z \in \tilde{Z}$. They have natural generalization to the case of secondary fan, where $Z \in \hat{Z}$. Such generalization turns out to be extremely straightforward. The only difference in the argument when $Z \in \hat{Z}$ is that for each $S = \{ m_0,m_1,m_2\} \in Z^{\rm top}$, $\{z_i = \delta^{w_{m_i}} a_{m_i} s_{m_i}\}_{i=0}^2$ defines an open covering (instead of embedding) $U^S_\epsilon \hookrightarrow {\mathbb{CP}^2}$, where $[z_0,z_1,z_2]$ is the homogeneous coordinate of ${\mathbb{CP}^2}$. Local models in section 5.1 can be pull back using the open covering maps in the same way as using the open embeddings in the case of $Z \in \tilde{Z}$. With this understanding, it is easy to check that all arguments in the case of $Z \in \tilde{Z}$ can easily be adopted to the case of $Z \in \hat{Z}$. We have\\
\begin{theorem}
\label{ef}
Theorems \ref{dd} and \ref{ee} are also true when $Z \in \hat{Z}$.
\hfill$\Box$\\
\end{theorem}
As we did at the end of section 4.1, we may classify the fibres $\mu_r := C_0 \cap \hat{F}^{-1}(r)$ of the map $\hat{F}: C_0 \rightarrow \Gamma$ for $r\in \Gamma$ in the general case. In general, when $Z \in \hat{Z}$, $\mu_r$ can be several points when $r$ is an end point of $\Gamma$. $\mu_r$ can be several circles when $r$ is a smooth point of $\Gamma$. $\mu_r$ can be an abelian multiple cover of the ``$\Theta$'' shaped graph in the torus at the right of figure \ref{fig9} when $r$ is a 3-valent vertex of $\Gamma$. (The graph illustrated at the left of figure \ref{fig9} can be viewed as an example of such, which is a $(\mathbb{Z}_5)^2$-cover of the ``$\Theta$'' shaped graph.) In the special case when $Z \in \tilde{Z}$, $\mu_r$ is a point when $r$ is an end point of $\Gamma$. $\mu_r$ is a circle when $r$ is a smooth point of $\Gamma$. $\mu_r$ is the ``$\Theta$'' shaped graph when $r$ is a 3-valent vertex of $\Gamma$.\\

{\bf Examples:} Using these theorems, the images of degree $d=5$ curves in ${\mathbb{CP}^2}$ under the moment maps as illustrated in figure \ref{fig2} can be perturbed to the following\\
\begin{center}
\begin{picture}(200,180)(-30,0)
\thicklines
\multiput(82,149)(36,0){1}{\line(2,1){18}}
\multiput(82,149)(36,0){1}{\line(-2,1){18}}
\multiput(82,128)(36,0){1}{\line(0,1){21}}
\multiput(46,128)(36,0){2}{\line(2,-1){18}}
\multiput(82,128)(36,0){2}{\line(-2,-1){18}}
\multiput(64,98)(36,0){2}{\line(0,1){21}}
\multiput(28,98)(36,0){3}{\line(2,-1){18}}
\multiput(64,98)(36,0){3}{\line(-2,-1){18}}
\multiput(46,68)(36,0){3}{\line(0,1){21}}
\multiput(10,68)(36,0){4}{\line(2,-1){18}}
\multiput(46,68)(36,0){4}{\line(-2,-1){18}}
\multiput(28,38)(36,0){4}{\line(0,1){21}}
\multiput(-8,38)(36,0){5}{\line(2,-1){18}}
\multiput(28,38)(36,0){5}{\line(-2,-1){18}}
\multiput(10,8)(36,0){5}{\line(0,1){21}}
\thinlines
\put(-26,8){\line(1,0){216}}
\put(-26,8){\line(3,5){108}}
\put(190,8){\line(-3,5){108}}
\end{picture}
\end{center}
\begin{center}
\stepcounter{figure}
Figure \thefigure: $F(C_p)$ of degree $d=5$ curve in ${\mathbb{CP}^2}$ perturbed to $\hat{F}(C_p) = \Gamma$
\end{center}
This example correspond to the large complex limit with respect to the standard simplicial decomposition of $\Delta$. When approaching different large complex limit in $\overline{{\cal M}_g}$ the toric moduli space of stable curves of genus $g$, the graph $\Gamma$ will be different and determined by the corresponding simplicial decomposition $Z$ of $\Delta$. Following is an example for degree $d=5$ curve in ${\mathbb{CP}^2}$.\\
\begin{center}
\begin{picture}(200,200)(-30,0)
\thicklines
\multiput(82,149)(36,0){1}{\line(2,1){18}}
\multiput(82,149)(36,0){1}{\line(-2,1){18}}
\multiput(82,149)(5,-31){2}{\line(2,-3){14}}
\multiput(118,129)(-31,-10){2}{\line(-1,0){23}}
\multiput(86,119)(-27,-10){1}{\line(1,1){10}}
\multiput(46,128)(36,0){1}{\line(2,-1){18}}
\multiput(82,128)(36,0){0}{\line(-2,-1){18}}
\multiput(64,98)(36,0){1}{\line(0,1){21}}
\multiput(28,98)(72,0){2}{\line(2,-1){18}}
\multiput(64,98)(27,-20){2}{\line(1,-2){9}}
\put(73,79){\line(1,0){18}}
\multiput(64,98)(72,0){2}{\line(-2,-1){18}}
\multiput(100,98)(-27,-20){2}{\line(-1,-2){9}}
\multiput(46,68)(72,0){2}{\line(0,1){21}}
\multiput(10,68)(36,0){2}{\line(2,-1){18}}
\multiput(118,68)(36,0){1}{\line(2,-1){18}}
\multiput(46,68)(36,0){1}{\line(-2,-1){18}}
\multiput(118,68)(36,0){2}{\line(-2,-1){18}}
\multiput(28,38)(36,0){4}{\line(0,1){21}}
\multiput(-8,38)(36,0){5}{\line(2,-1){18}}
\multiput(28,38)(36,0){5}{\line(-2,-1){18}}
\multiput(10,8)(36,0){5}{\line(0,1){21}}
\thinlines
\put(-26,8){\line(1,0){216}}
\put(-26,8){\line(3,5){108}}
\put(190,8){\line(-3,5){108}}
\end{picture}
\end{center}
\begin{center}
\stepcounter{figure}
Figure \thefigure: Alternative $\Gamma$ for degree $d=5$ curve in ${\mathbb{CP}^2}$
\end{center}
Applying these theorems to the case of curves in the toric surface (${\mathbb{CP}^2}$ with 3 points blown up) as illustrated in figure \ref{fig3}, we will be able to perturb the image of the moment map to the following graph.\\
\begin{center}
\begin{picture}(200,130)(-30,0)
\put(82,128){\makebox(0,0){$E_1$}}
\put(162,23){\makebox(0,0){$E_3$}}
\put(4,23){\makebox(0,0){$E_2$}}
\thicklines
\multiput(64,98)(36,0){2}{\line(0,1){21}}
\multiput(28,98)(36,0){3}{\line(2,-1){18}}
\multiput(64,98)(36,0){3}{\line(-2,-1){18}}
\multiput(46,68)(36,0){3}{\line(0,1){21}}
\multiput(10,68)(36,0){4}{\line(2,-1){18}}
\multiput(46,68)(36,0){4}{\line(-2,-1){18}}
\multiput(28,38)(36,0){4}{\line(0,1){21}}
\multiput(28,38)(36,0){4}{\line(2,-1){18}}
\multiput(28,38)(36,0){4}{\line(-2,-1){18}}
\multiput(46,8)(36,0){3}{\line(0,1){21}}
\thinlines
\put(40.6,119){\line(1,0){82.8}}
\put(22.25,8){\line(1,0){119.5}}
\put(22.25,8){\line(-3,5){24}}
\put(141.75,8){\line(3,5){24}}
\put(40.6,119){\line(-3,-5){42.5}}
\put(123.4,119){\line(3,-5){42.5}}
\end{picture}
\end{center}
\begin{center}
\stepcounter{figure}
Figure \thefigure: $F_s(C_s)$ in figure \ref{fig3} perturbed to graph $\Gamma$
\vspace{10pt}
\end{center}
{\bf Acknowledgement:} I would like to thank Prof. S.T. Yau for constant encouragement, Prof. Yong-Geun Oh for pointing out the work of \cite{M} to me. This work was initially done while I was in Columbia University. I am very grateful to Columbia University for excellent research environment. Thanks also go to Qin Jing for stimulating discussions and suggestions.\\\\

\end{document}